%% file: papersl3R.tex
\documentclass{amsart}
\usepackage{amsmath,amsthm,amsfonts,amssymb, mathrsfs}
\usepackage{graphicx,xcolor}
\usepackage{comment}
\usepackage{todonotes}
\usepackage{tikz-cd}
\usepackage[all]{xy}
\usepackage[]{todonotes}
\usepackage{svg}

\usepackage[backend=biber,style=alphabetic,safeinputenc]{biblatex}
\newcommand{\RR}{\mathbb{R}}

\newcommand{\HH}{\mathbb{H}}
\newcommand{\PP}{\mathbb{P}}

\DeclareMathOperator{\trace}{tr}
\DeclareMathOperator{\tess}{Tess}
\DeclareMathOperator{\dev}{Dev}

\newtheorem{thm}{Theorem}[subsection]
\newtheorem{lemma}[thm]{Lemma}
\newtheorem{prop}[thm]{Proposition}
\newtheorem{coro}[thm]{Corollary}

\theoremstyle{remark}
\newtheorem{rmk}[thm]{Remark}

\theoremstyle{definition}
\newtheorem{defi}[thm]{Definition}

\newtheorem{thmintro}{Theorem}

\title{Parallel domains in $SL_3(\mathbb{R})$ and the ideal triangle reflection group}

\author{Joaqu\'in Lejtreger}
\address[Joaqu\'in Lejtreger]{
Max Planck Institute For Mathematics in the Sciences, Inselstra{\ss}e 22, 04103 Leipzig}
\email{joaquin.lejtreger@mis.mpg.de}
\date{}

\author{Joaqu\'in Lema}
\address[Joaqu\'in Lema]{Department of Mathematics, Boston College, Chestnut Hill, MA 02467}
\email{lemajo@bc.edu}
\date{}

\begin{document}

\begin{abstract}
For a given positive triple of flags, we associate an open set, called a parallel domain, in the symmetric space of $SL_3 (\mathbb{R})$ whose boundary is parametrized by three spacelike geodesics of the pseudo-Riemannian symmetric space $SL_3 (\mathbb{R})/GL_2 (\mathbb{R})$. This construction is a higher-rank analog of how three points in the ideal boundary of the hyperbolic plane have an associated ideal triangle.

As an application, using the fact that positive triples determine a representation of the ideal triangle reflection group into $SL_3(\mathbb{R})$, we show that these parallel domains are fundamental domains for the group action in the symmetric space. 
Moreover, using the explicit construction and combinatorial properties of the Farey tessellation, we compute the asymptotics for the critical exponent in the symmetric space and the Hilbert critical exponent as the triple ratio of the triple of flags goes to infinity.
\end{abstract}

\maketitle

\input{sections/intro}     
\input{sections/symmetric} 
\input{sections/walls}
\input{sections/triples.tex}
\input{sections/domains}
\input{sections/criticalexponent.tex}

\appendix

\input{sections/negative.tex}

\printbibliography

\end{document}

%% file: sections/intro.tex
\section*{Introduction}

Any three distinct points on the boundary of the hyperbolic plane span a unique ideal triangle, and any two ideal triangles are isometric.
Gluing these rigid pieces with a shear on each side gives coordinates on the Teichm\"uller space of finite-area hyperbolic surfaces \cite{bonahon1996shearing}, and even on the universal Teichm\"uller space \cite{saric2010circle}.

This picture has been influential in higher rank.
Fock--Goncharov parametrized convex projective structures on punctured surfaces by assigning positive configurations of flags to the punctures \cite{fock2007moduli}, and these coordinates were later extended to closed surfaces and to Hitchin components of higher rank \cite{fock2004moduli,bonahon2014parameterizing}.
While these coordinates help study Hitchin components, their connection to the geometry of the symmetric space is unclear.

In this paper, we define an analog of the ideal triangle for $G = SL_3(\RR)$.
Given a positive triple of flags $(F_1, F_2, F_3)$, we associate a domain $\mathcal{D}(F_1, F_2, F_3)$ in the symmetric space $\mathbb{X}$, bounded by three hypersurfaces that are asymptotic to the three flags, in an appropriate sense.
When the flags are tangent to an ellipse in $\RR\PP^2$, the domain coincides with the preimage of an ideal triangle under the closest point projection onto a totally geodesic hyperbolic plane.

The main difficulty in finding such a domain is that the symmetric space has no totally geodesic hypersurfaces, so it has no natural analog of the sides of an ideal triangle.
Our hypersurfaces, which we call \emph{walls}, are thus not totally geodesic, but they are foliated by totally geodesic submanifolds called \emph{parallel sets}.
This is reminiscent of bisectors in complex hyperbolic geometry, which are foliated by complex geodesics \cite{goldman1999complex}.

As an application of the domains, we study convex projective structures on the orbifold associated to the group generated by reflections along the sides of an ideal triangle of $\HH^2$, which we denote by $W_3$.
The holonomy of these convex projective structures coincides with the relative Hitchin representations of $W_3$ into $SL_3(\RR)$ (see \cite{choi2005deformation} in the compact setup and \cite{marquis2010espace} in the non-compact setup). Positive triples of flags parametrize these representations.

\subsection*{The setup}

An irreducible hyperbolic slice $\mathcal{H}$ in the symmetric space $\mathbb{X}$ is a totally geodesic copy of $\HH^2$ whose stabilizer in $SL_3(\RR)$ is $SO(B)$, for $B$ some $(1,2)$ form in $\RR^3$.
A key observation is that the closest-point projection $\pi: \mathbb{X} \to \mathcal{H}$ defines a fibration of $\mathbb{X}$, with totally geodesic fibers called \emph{parallel sets}.

Parallel sets are defined as the union of flats containing a given singular geodesic. Each parallel set is isometric to $\HH \times \RR$, where $\HH$ is a totally geodesic copy of the hyperbolic plane stabilized by a reducible embedding of $SL_2(\RR)$ into $SL_3(\RR)$. The space of parallel sets is parametrized by the homogeneous space
$$\mathbb{Y} = \{(p, \ell) \in \RR\PP^2 \times (\RR\PP^2)^*: p \notin \ell\}.$$
This space has plenty of structure; for instance, it possesses a natural symplectic form as described in \cite{labourie2007cross}. 
We will focus on its structure as a pseudo-Riemannian symmetric space. 
Crucially, two parallel sets joined by a spacelike geodesic are disjoint and their distance in $\mathbb{X}$ can be computed using the cross ratio (see Proposition \ref{prop:relativeposition}).

\begin{figure}[ht]
    \centering
    \includegraphics[scale=1]{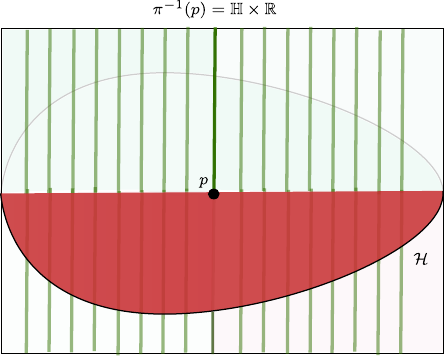}
    \caption{Schematic of a wall. The red region is an irreducible slice, and the wall divides $\mathbb{X}$ into two regions.}
    \label{fig:walls}
\end{figure}

Given a geodesic $\gamma$ in $\mathcal{H}$, $\pi^{-1}(\gamma)$ is the union of the parallel sets corresponding to a spacelike geodesic in $\mathbb{Y}$. Conversely, any spacelike geodesic in $\mathbb{Y}$ is of this form. 
The union of these parallel sets is a $4$-dimensional hypersurface of $\mathbb{X}$ which we call a \emph{wall}.

Spacelike geodesics in $\mathbb{Y}$ are asymptotic to a pair of flags in general position. Conversely, given a pair of flags in general position, there are two spacelike geodesics asymptotic to them (see Proposition \ref{prop:asymptoticgeodesics}).

\begin{figure}[ht]
    \centering
    \includegraphics[scale=.8]{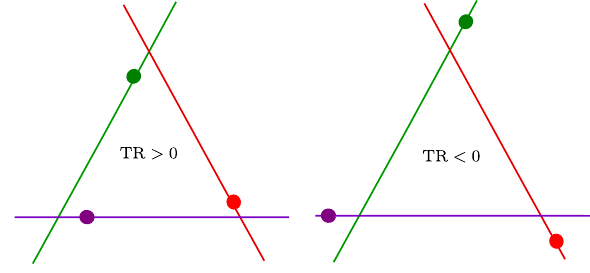}
    \caption{The triple ratio is a projective invariant of a triple of flags. The figure depicts a positive and a negative triplet.}
    \label{fig:triplet-tr}
\end{figure}

Triples of flags in general position come in two kinds, positive and negative (see Figure \ref{fig:triplet-tr}).
We prove that, given a positive triple of flags, there is a canonical choice of spacelike geodesics for each pair. Moreover, we show that the walls associated to these three geodesics are disjoint and bound an open set in $\mathbb{X}$ which we call the \emph{parallel domain} associated to the triple.
When the three flags are tangent to an ellipse, the parallel domain is the preimage of an ideal triangle under the closest-point projection onto an irreducible hyperbolic slice.

A useful feature of the parallel domain is that, using cross ratios, we can explicitly compute distances between parallel sets contained in walls in the boundary of the domain (see Proposition \ref{prop:estimates}). 

\subsection*{Main results}

We apply the ideas explained above to study representations of the group $W_3$. A presentation for this group is
$$W_3 = \left\langle s_1, s_2, s_3 \mid s_i^2 \right\rangle.$$
More specifically, we study representations of this group into $SL_3(\RR)$ in which the products of pairs of distinct generators are unipotent and fix a unique flag. 
These representations are parametrized by triples of flags fixed by these three products. 
With this in mind, we prove the following theorem:

\begin{thmintro}[Theorem \ref{thm:main}]
\label{thmintro1}
    Let $(F_1, F_2, F_3)$ be a positive triple of flags, and let $\rho: W_3 \to SL_3(\RR)$ be its associated representation. Then the parallel domain $\mathcal{D}(F_1, F_2, F_3)$ is a fundamental domain for the action of $W_3$ on $\mathbb{X}$ via $\rho$.
\end{thmintro}

Notice that discreteness and faithfulness is not new, as this follows from classical work of Vinberg \cite{vinberg1971discrete}.
The novelty is the explicit construction of a fundamental domain in the symmetric space.

The proof follows the classical strategy used to prove the Poincar\'e polyhedron theorem, or Vinberg's criterion (see \cite{epstein1994exposition,benoist2009five}).
One first builds a tesselated space on which $W_3$ acts with $\mathcal{D}$ as a fundamental domain, and equipped with a local isometry to the symmetric space.
Then, one shows that the tesselated space is complete.
This step requires the estimates with the parallel domains that we mentioned before. 

Our treatment relies solely on the geometry of the symmetric space.
One can also prove Theorem \ref{thmintro1} using results from the theory of geometrically finite Hilbert geometries \cite{crampon2014finitude}.
A proof along those lines can be found in the first author's thesis \cite{lejtreger2026geometric}.

Orbits of triples of flags in general position under the diagonal action of $SL_3(\RR)$ are parametrized by the \emph{triple ratio}. 
The relative character variety sending the peripheral system of the reflection orbifold to regular unipotents $\mathcal{X}_{\text{rel}}^{\text{irr}}(W_3, SL_3(\RR))$ is thus parametrized by $\mathbb{R}\setminus \{0, -1\}$ via the triple ratio. 
Positive representations correspond to positive triple ratios.

We then proceed to use these estimates to compute the asymptotics of the critical exponent of $\rho_\tau$ associated to a triple of flags with triple ratio $\tau$ in the symmetric space.
We define the \emph{critical exponent of $\rho$} in the symmetric space to be
$$h_\mathbb{X} (\rho) = \limsup_{T \to \infty} \frac{1}{T}\log \#\{g \in W_3: d_\mathbb{X}(o, \rho(g)o) \le T\},$$
for some choice of basepoint $o \in \mathbb{X}$. 
The number is independent of this choice.

\begin{thmintro}[Theorem \ref{thm:limite-exp}]
\label{thm:critintro}
For each $\tau > 0$, we let $\rho_\tau$ be a representation given by a triple of flags with triple ratio $\tau$. Then
$$\lim_{\tau \to +\infty} h_\mathbb{X} (\rho_\tau) = \frac{1}{2}.$$
\end{thmintro}

We prove the (stronger) statement that the Hilbert critical exponent converges to $\frac{1}{2}$ (see Theorem \ref{thm:criticalhilbert}). 
Giuseppe Martone, Harry Bray and Mitul Islam supervised an undergraduate project \cite{o2021degeneration} in which the authors estimated the Hilbert critical exponent numerically.
In that article, they conjecture the $1/2$ limit at infinity, which the previous result confirms.

For $\tau = 1$, the representation preserves an irreducible slice $\mathcal{H}$ with curvature $-1$, so the critical exponent is one.
Theorem \ref{thm:limite-exp} says that as $\tau$ goes to infinity, the critical exponent drops to $1/2$, which coincides with the critical exponent of a cyclic subgroup generated by a regular unipotent element.
There are two key ingredients in the proof of the theorem above.
The first is to use distance estimates in the parallel domain to produce lower bounds for the translation length of elements with little excursion along the cusps compared to their word length.
The second ingredient is to obtain a bound for the translation length of those elements with large excursions on the cusps, which we reduce to a computation in the hyperbolic plane (see Section \ref{subsec:farey}). Moreover, our computations hint to a positive answer to the following question (see Proposition \ref{prop:monotone} and Theorem \ref{thm:coxeteriality}).

\smallskip

\textbf{Question:} Is the function $\tau \to h_\mathbb{X} (\rho_\tau)$ strictly decreasing in $[\tau_0,\infty)$ for some positive real $\tau_0$?

\smallskip

One can observe that there is a duality that exchanges $\rho_\tau$ with $\rho_{1/\tau}$. Thus, analizing the case $\tau \to \infty$ is analogous to $\tau \to 0$.

\subsection*{Related work}
Our motivation came from Schwartz's work on Pappus representations \cite{schwartz2025leretourpappus}, where he associates a pattern of flats with a Pappus representation and builds a sort of pleated surface from it.
In our setup, we arrived at the construction of domains by examining configurations of flats foliated by medial geodesics reminiscent of his work. 
We are able to fatten these patterns to a hypersurface due to the positivity of the flag configurations of the representations studied.

There is an extensive list of works providing asymptotics for critical exponents in the Hitchin component, dating back to Nie \cite{nie2015entropy} and Zhang \cite{zhang2015degeneration}.
Recently, in \cite{carvajales2026triangular} the authors provide asymptotics for specific families of Hitchin representations in terms of Fock--Goncharov parameters on a compact surface with boundary, based on work by Parreau \cite{parreau2022invariant}.
The setup of the previous work is different from ours as they avoid working in the symmetric space, relying on projective techniques instead.
In previous work \cite{lejtreger2026separated} we study entropy asymptotics for convex projective structure on pairs of pants, but along sequence with bounded triple ratio, contrary to the case studied here.
In upcoming work, we will extend our domains in order to tackle the compact case, by providing a reinterpretation of the Bonahon--Dreyer coordinates in $SL_3 (\mathbb{R})$ \cite{bonahon2014parameterizing}.

Recently, the pseudo-Riemannian structure of $\mathbb{Y}$ has been exploited in \cite{wang2026geometric} in relation to the Axiom A flows constructed in \cite{delarue2025locally}.
\subsection*{Organization of the paper}

In Section \ref{domains:sec:symmetric} we review classical theory on the symmetric space of $SL_3(\RR)$.
In Section \ref{domains:sec:walls}, we explore the pseudo-Riemannian structure of the space $\mathbb{Y}$ and introduce the notion of a wall.

Section \ref{sec:parallel} is the technical heart of the paper. 
In it, we construct the parallel domains, compute distances between parallel sets in the walls on the boundary, and obtain other important estimates. 

Using these results, Section \ref{domains:sec:representations} proves that the parallel domains are fundamental domains for the $W_3$ representations we consider.
Section \ref{sec:critical-exponent} estimates the critical exponents as the triple ratio goes to infinity.

We include Appendix \ref{sec:negative} with some questions about the negative component of the relative character variety of $W_3$. 
We show that some of these representations are potential counterexamples to many currently open questions on discrete subgroups of $SL_3 (\mathbb{R})$.

\subsection*{Acknowledgments}
We thank Rich Schwartz for inspiring conversations about symmetric spaces, Giuseppe Martone and Sara Maloni for their interest in our work and for sharing theirs, and François Gu\'eritaud and Andrea Monti for their useful comments on an earlier version of this paper. 
We are grateful to both Andr\'es Sambarino and Martin Bridgeman for their guidance and support, and to Le\'on Carvajales for interesting conversations through the years.

\subsection*{Tool use Statement}
We used GenAI to proofread a preliminary version of this article, which detected typos and gaps in proofs that we then corrected manually.
All mathematical arguments were developed without AI assistance.
We created most figures using only Inkscape, except Figures \ref{fig:farey}, \ref{fig:negative}, which we created with assistance from AI-generated code.

%% file: sections/symmetric.tex
\section{The symmetric space of $SL_3(\RR)$}
\label{domains:sec:symmetric}

In this section we will study the geometry of the symmetric space of $SL_3(\RR)$, which we denote by $\mathbb{X}$.
We will only take a selective overview of the topic. 
The reader can find a more thorough treatment for the $SL_n (\mathbb{R})$-case in \cite[Section 2]{schwartz2025leretourpappus} or \cite[Chapter 4]{donaldson2007lie}, and \cite{sambarino2024espacios}, \cite{helgason1979differential}, \cite{eberlein1996geometry} for generalities on symmetric spaces.

We define $\mathbb{X}$ as the space of unimodular inner products on $\mathbb{R}^3$, i.e., inner products whose unit ball has the same volume as the standard Euclidean ball. 
The group $SL_3 (\RR)$ acts transitively on $\mathbb{X}$ via pullback $g \cdot q = (g^{-1})^* q$, with stabilizer $SO(q)$ (by definition). 

\subsection{The metric}
\label{subsec: generalities}

Denote by $\mathfrak{p}_q$ the space of $q$-self adjoint, traceless matrices. The assignment 
$$X \mapsto \frac{d}{dt}\big|_{t = 0} e^{tX}\cdot q$$
identifies $T_q \mathbb{X} \cong \mathfrak{p}_q$. 
Note that the derivative of the action by $g$ is $\mathrm{Ad}(g)$.

We can equip $\mathbb{X}$ with a Riemannian metric by declaring
$$\langle X,Y \rangle_q = \frac{1}{2}\trace (X Y),\: \forall X,Y \in \mathfrak{p}_q.$$
Because of the description of the derivative of the action of $SL_3(\RR)$, we get that $SL_3(\RR)$ acts by isometries. 
More specifically, the metric identifies the group with $\mathrm{Isom}_0(\mathbb{X})$. 
The metric coincides, up to some multiple, with the metric inherited from the Killing form of $SL_3 (\RR)$.
Our choice of multiple is explained by Remark \ref{rmk:curvature}.

One can prove that the geodesic through $q$ tangent to $X \in \mathfrak{p}_q$ is
$$\gamma_{(q,X)} (t) = e^{tX} \cdot q.$$

Geodesics through $q \in \mathbb{X}$ come in two flavors, depending on the velocity vector $X \in \mathfrak{p}_q$.
Since $X$ is $q$-self-adjoint, it is diagonalizable in a $q$-orthonormal basis.
If $X$ has three distinct eigenvalues, we say the geodesic is \emph{regular}.
Otherwise, we say that the geodesic is \emph{singular}.

\subsection{Curvature}
If we denote by $K$ the sectional curvature, one can prove the following formula for the sectional curvature:

\begin{equation}
\label{eq:curvature}
K (X, Y) = \frac{2\trace ([X,Y]^2)}{\trace(X^2) \trace(Y^2) -  \trace(XY)^2}.
\end{equation}

A quick computation shows that this number is always non-positive, and in particular, $\mathbb{X}$ is a $CAT(0)$-space.

\subsection{Regular geodesics, flats and flags}
\label{subsec:regular}
Given $q \in \mathbb{X}$ and a $q$-orthonormal basis $\mathcal{B}$, we can consider a subspace $\mathfrak{a}_\mathcal{B} \subseteq \mathfrak{p}_q$ of matrices diagonal in $\mathcal{B}$. 
One can verify that the image of the map
$$
\begin{aligned}
 \mathfrak{a}_\mathcal{B} &\rightarrow \mathbb{X}\\
 z &\mapsto e^z \cdot q
\end{aligned}
$$ 
is a totally geodesic surface isometric to the euclidean plane $\mathbb{E}^2$ (up to some multiple). We will call such a plane a \textit{flat}. 

The action of $SO(q)$ is transitive on the set of flats through $q$. 
A flat $A$ is then determined by three points on $\RR\PP^2$ in general position, since such a triple defines an orthogonal basis for each inner product in $A$.
If $X$ is tangent to a regular geodesic, there is a unique $q$-orthogonal basis in which $X$ is diagonal. 
Hence, it lies in a unique flat. 
A singular geodesic is contained in a circle worth of flats. 

When the geodesic is regular, the eigenspaces have a natural order, namely, which eigenspace corresponds to the larger eigenvalue.
This information can be codified with a pair of \emph{flags} in general position.

A flag is a pair $(p,\ell)$ consisting of a one-dimensional subspace $p \in \mathbb{RP}^2$, contained in a two-dimensional subspace $\ell \in (\mathbb{RP}^2)^*$.
We will denote the space of flags as $\mathcal{F}$.
A pair of flags $F_1, F_2 \in \mathcal{F}^2$ are said to be in \emph{general position}, or \emph{transverse} if 
$$p_1 + p_2 + \ell_1 \cap \ell_2 = \mathbb{R}^3.$$
We denote the space of pairs of flags in general position as $\mathcal{F}^{(2)}$.

A regular $X \in \mathfrak{p}_q$ determines a pair of flags $F_X^+$ and $F_X^-$, whose one-dimensional subspaces correspond to the top (resp. bottom) eigenspace, and whose two-dimensional subspace contains both the top (resp. bottom) and middle eigenlines.
Such a pair $(F_X^+,F_X^-)$ is in general position by construction.

Since $\mathbb{X}$ is a non-positively curved, simply connected manifold, we can construct its visual boundary $\partial \mathbb{X}$.
Inside the visual boundary, we get $\partial^{\text{reg}} \mathbb{X} \subset \partial \mathbb{X}$, the points at infinity defined by regular geodesic rays.
The construction of the flags $F_X^+$ from before defines a map $\partial^{\text{reg}} \mathbb{X} \rightarrow \mathcal{F}$ by sending the point corresponding to a geodesic ray $[\gamma_X]$ associated to $X \in \mathfrak{p}_{\gamma_X (0)}$ to $F_X^+$.
The preimage of a flag $F$ consists of a copy of:

$$\mathfrak{a}^{++} = \{(x,y,z) \in \mathbb{R}^3: x+y +z = 0,\: x > y > z\},$$

called the \emph{interior of a standard Weyl chamber}.
The correspondence is given by taking all $X \in T \mathbb{X}$ that are regular and have $F_X^+ = F$, and collecting the logs of their top and middle eigenvalues.
The unit-determinant condition automatically determines the smallest eigenvalue.
More generally, the \emph{Weyl-chamber} $\mathfrak{a}^+$ is the closure of $\mathfrak{a}^{++}$.

We will refer to the subspace $\mathfrak{a} \subset \mathbb{R}^3$ defined by $x+y + z =0$ as the \emph{Cartan subspace}.
We will equip this subspace with the metric 
$$\| (x,y,z)\|^2 = \frac{1}{2} (x^2 + y^2 + z^2).$$
With this normalization, given $q_0$ the standard inner product in $\mathbb{R}^3$, we can embed $\mathfrak{a} \rightarrow \mathfrak{p}_{q_0}$ as $(x,y,z) \to \text{diag}(x,y,z)$, and this map is isometric onto its image with our choice of metric.
In fact, the exponential of the image defines a flat passing through $q_0$ which we call \emph{the standard flat}.

Reflections along the faces of $\mathfrak{a}^+ \subset \mathfrak{a}$ generate the symmetric group $S_3$ (the \emph{Weyl} group).
Thinking of $\mathfrak{a}$ inside $\mathfrak{p}_{q_0}$, one can prove that these reflections are generated by the differentials of elements that stabilize $q_0$.

\subsection{Cartan projections}
\label{subsec:finslercartan}

We start by describing a \emph{vector-valued distance} in the symmetric space $\vec{d} : \mathbb{X} \times \mathbb{X} \rightarrow \mathfrak{a}^+$.
Given $q,p \in \mathbb{X}$, we can always write $p = \exp(X).q$, where $X \in \mathfrak{p}_q$ is the vector defining the geodesic joining the two points.
As before, the vector $X$ is diagonalizable with eigenvalues $\lambda_1 \geq \lambda_2 \geq \lambda_3$.
We define $\vec{d} (q,p)$ to be the element $(\lambda_1,\lambda_2,\lambda_3) \in \mathfrak{a}^+$.

Geometrically, we can always find an isometry that maps the pair $(q,p)$ to $(q_0,p_0)$, where $q_0$ is the standard inner product and $p_0$ is in the standard flat.
The vector-valued distance corresponds to the vector corresponding to $p_0$ after identifying the flat with $\mathfrak{a}$, and maybe moving $p_0$ through the stabilizer of $q_0$ until it lies in the image of $\mathfrak{a}^+$.

Fixing a basepoint $q \in \mathbb{X}$, we define the induced \emph{Cartan projection} $\kappa : SL_3 (\mathbb{R}) \rightarrow \mathfrak{a}^+,$ as $\kappa(g) =\vec{d}(q,gq)$.

\subsection{Dualities and irreducible slices}
\label{subsub:hypslices}

Let $B$ be a non-degenerate symmetric bilinear form (of arbitrary signature).
Such a tensor can be used to define a map of $\mathbb{X}$ via:
$$I_B (q) (v,w) = q^* (B(v,\cdot),B(w,\cdot)) = q([B]_q\cdot,[B]_q\cdot),$$
where $q^*$ is the induced inner product in $(\mathbb{R}^3)^*$, and $[B]_q$ is the element such that $B(\cdot,\cdot) = q ([B]_q \cdot,\cdot)$.
This is well-defined as long as $\det [B]_q= 1$ for some (thus, for all) $q \in \mathbb{X}$ which we assume from now on.
The determinant condition forces $B$ to be an inner product, or have signature $(1,2)$.

Notice that for any $g \in G$, we get that $I_B (g q) = (g^{*_B})^{-1} I_B(q)$, for $*_B$ the adjoint with respect to $B$.
In particular, we get that $I_B$ is an involution, and $d_q I_B (X) = - X^{*_B}$.
This second condition implies that $I_B$ is an isometry.
Thus, its fixed point set is totally geodesic.

Notice that when $q \in \mathbb{X}$ is a point in the symmetric space, the map $I_q$ is the \emph{symmetric involution} around $q$, that is, the map fixing $q$ and acting as $-Id$ in $T_q \mathbb{X}$.

\begin{rmk}
\label{rmk:compinvolutions}
Given $q \in \mathbb{X}$ and $B$ a non-degenerate bilinear forms with well-defined $I_B$, then $I_{B} I_{q} I_{B} = I_{I_B (q)}$.
In particular, if $I_B (q) = q$ we get that $I_B I_q$ is again, an involution.
Moreover, $I_q I_B$ equals multiplication by $[B]_q$.
\end{rmk}

When $B$ is not positive definite, we denote the fixed point set of $I_B$ by $\mathcal{H}_B$.
More concretely, we can write:
$$\mathcal{H}_B = \{q \in \mathbb{X}: [B]_q \in SO(q)\}.$$
We will call such a submanifold an \emph{irreducible hyperbolic slice}.
Being a totally geodesic submanifold, it is necessarily a symmetric space.
Moreover, $SO(B)$ acts on the fixed point set $\mathcal{H}_B$, and it acts transitively.
In particular, $\mathcal{H}_B$ must be isometric to the symmetric space of $SO(B)$, which is the hyperbolic plane.

\begin{rmk}
\label{rmk:curvature}
Applying Equation \ref{eq:curvature} shows that, with our normalizations, the slice has constant curvature $-1$.
\end{rmk}

%% file: sections/walls.tex
\section{Parallel sets and walls in the symmetric space}
\label{domains:sec:walls}

In this section we study the space $\mathbb{Y}$ of \emph{parallel sets}, a special family of totally geodesic submanifolds of $\mathbb{X}$.
Each parallel set is the fixed point set of an involution of $SL_3(\RR)$, and this identifies $\mathbb{Y}$ with a conjugacy class of involutions.
In particular, $\mathbb{Y}$ is a pseudo-Riemannian symmetric space.
Its spacelike geodesics correspond to hypersurfaces of $\mathbb{X}$ saturated by parallel sets, which we call \emph{walls}.
 
\subsection{Linear involutions and the symmetric space structure}
\label{prelim:involutions}

A \emph{linear involution} of a real vector space $V$ is an element $\sigma \in \operatorname{GL}(V)$ such that $\sigma^2 = \mathrm{Id}$.
We denote by $\operatorname{Gr}_k(V)$ the Grassmannian of $k$-dimensional subspaces of $V$, and write $\mathbb{P}(V) = \operatorname{Gr}_1(V)$ and $\mathbb{P}(V)^* = \operatorname{Gr}_{\dim V - 1}(V)$.
 
\begin{lemma}
\label{lemma:involutions}
Every linear involution $\sigma$ of $V$ is diagonalizable with eigenvalues $\pm 1$.
Denoting by $\sigma^\pm$ the corresponding eigenspaces, the map $\sigma \mapsto (\sigma^+,\sigma^-)$ is a $\operatorname{GL}(V)$-equivariant bijection between the conjugacy class of a given involution $\mathcal{C}(\sigma_0) = \{ g \sigma_0 g^{-1} : g \in \operatorname{GL}(V)\}$ and
$$\mathbb{Y}_\sigma = \{ (V^+,V^-) \in \operatorname{Gr}_{\dim \sigma^+} (V) \times \operatorname{Gr}_{\dim \sigma^-} (V) : V^+ \cap V^- = 0\}.$$
\end{lemma}
 
\begin{proof}
Every $v \in V$ decomposes as $v = v^+ + v^-$, where $v^{\pm} = \frac{1}{2}(v \pm \sigma v)$ satisfies $\sigma v^{\pm} = \pm v^{\pm}$.
This proves the first claim.
 
For $g \in \operatorname{GL}(V)$, the eigenspaces of $g\sigma g^{-1}$ are $g\sigma^{\pm}$, so the map is equivariant.
The map $\sigma \to (\sigma^+,\sigma^-)$ is clearly injective and it is surjective since $\operatorname{GL}(V)$ acts transitively on pairs of complementary subspaces of fixed dimensions.
\end{proof}
 
From now on $V = \RR^3$ and $G = SL_3(\RR)$, with Lie algebra $\mathfrak{g} = \mathfrak{sl}_3(\RR)$.
If $\sigma \in G$ is an involution then the unit determinant condition forces $\dim \sigma^+ = 1$ and $\dim \sigma^- = 2$.
We set
$$\mathbb{Y} := \{ (p,\ell) \in \RR\PP^2 \times (\RR\PP^2)^* : p \notin \ell\}.$$
By Lemma \ref{lemma:involutions}, $\mathbb{Y}$ is in $G$-equivariant bijection with the set of involutions of $G$.
We denote by $\sigma_{(p,\ell)}$ the involution associated to $(p,\ell) \in \mathbb{Y}$.
 
\begin{rmk}
\label{rmk:explicitinvolution}
The involution $\sigma_{(p,\ell)}$ can be written explicitly.
Pick $v \in p$ and $\alpha \in (\RR^3)^*$ with $\ker \alpha = \ell$ and $\alpha(v) = 2$, then
$$\sigma_{(p,\ell)} (w) = \alpha(w)v - w.$$
The pair $(v,\alpha)$ is unique up to replacing it by $(\lambda v, \lambda^{-1}\alpha)$ for $\lambda \in \RR \setminus \{0\}$.
\end{rmk}

\subsubsection*{Tangent space and the pseudo-Riemannian structure}
Using the splitting $\RR^3 = p \oplus \ell$, the tangent space of $\mathbb{Y}$ can be identified with
$$T_{(p,\ell)} \mathbb{Y} = \operatorname{Hom}(p,\ell) \oplus \operatorname{Hom} (\ell,p).$$
We can think of $T_{(p,\ell)} \mathbb{Y}$ inside $\operatorname{End} (\mathbb{R}^3)$ by using the splitting
$$\operatorname{End} (\mathbb{R}^3) = \left( \operatorname{End} (p) \oplus \operatorname{End}(\ell) \right) \oplus \left(  \operatorname{Hom}(p,\ell) \oplus \operatorname{Hom} (\ell,p)\right),$$
where the first term in the direct sum is fixed by $Ad (\sigma_{(p,\ell)})$, whereas the second term is anti-fixed by it.

Equipping $\mathfrak{g}$ with a multiple of the Killing form
$$\langle X, Y \rangle = \frac{1}{2}\operatorname{tr}(XY),$$
for any $X, Y \in \mathfrak{g}$, gives a pseudo-Riemannian metric on $\mathbb{Y}$ via the decomposition from before.
In particular, given $(S,T) \in T_{(p,\ell)} \mathbb{Y}$, we get
$$\langle (S,T),(S,T)\rangle = \text{tr}(T \circ S),$$

It is easy to see that this bilinear form has signature $(2,2)$ and the $G$ action is isometric by construction.

\subsubsection*{Symmetric involutions and geodesics}

Notice that $\sigma_{(p,\ell)}$ acts on $\mathbb{Y}$ fixing $(p,\ell)$.
Moreover, the differential of the action at $T_{(p,\ell)}\mathbb{Y}$ is precisely $\operatorname{Ad}(\sigma_{(p,\ell)})$ which, as we saw before, is $-Id$.
For that reason, we will call $\sigma_{(p,\ell)}$ the \emph{symmetric involution} of $\mathbb{Y}$ corresponding to $(p,\ell)$.

Given $\gamma$ a geodesic in $\mathbb{Y}$, a classical argument shows that $\sigma_{\gamma (t/2)} \sigma_{\gamma (0)}$ is a one parameter subgroup of the form $\exp (t X)$.
In particular, $\gamma (t) = \exp (tX) (p,\ell)$, for $X \in T_{(p,\ell)} \mathbb{Y}$.
We will now classify geodesics up to the $G$-action.

Notice that the $G$-stabilizer $H_{(p,\ell)}$ of a point $(p,\ell)$ is identified with $GL(\ell)$ thought of as linear transformations preserving the splitting $\mathbb{R}^3 = p \oplus \ell$.
Thus, to classify geodesics, it suffices to classify $H_{(p,\ell)}$ orbits of vectors in $T_{(p,\ell)} \mathbb{Y}$.

Let $\{e_i\}_{i=1,2,3}$ be the canonical basis of $\mathbb{R}^3$.
Using the homogeneous space structure, it is enough to classify the point $([e_1],[e_1^*]) \in \mathbb{Y}$.
Notice that in this case
$$H := H_{([e_1],[e_1^*])} = \left\{ 
    \begin{pmatrix}
        \frac{1}{\det g} & 0\\

        0 & g
    \end{pmatrix}: g \in GL_2 (\mathbb{R}) \right\} \cong GL_2 (\mathbb{R}).$$

Moreover, our choice of basepoint allows us to think of $T_{[e_1],[e_1^*]} \mathbb{Y}$ as $\langle e_2,e_3\rangle \oplus \langle e_2^*,e_3^* \rangle$.
The group $H$ acts through the following action: 
$$g \cdot (v,\alpha) = ((\det g) gv, \frac{1}{\det g} \alpha \circ g^{-1}).$$
A direct linear algebra computation shows the following lemma:

\begin{lemma}
\label{lemma:geodesicclassification}
Let $(v,\alpha) \in T_{([e_1],[e_1^*])} \mathbb{Y}$ be non-zero. 
Then $\langle (v,\alpha), (v,\alpha)\rangle = \alpha (v)$, and
\begin{itemize}
    \item if $(v,\alpha)$ is spacelike, then it is in the $H$-orbit of $\sqrt{\alpha (v)}(e_2,e_2^*)$.
    \item If $(v,\alpha)$ is lightlike, then it is in the $H$-orbit of exactly one of $(e_2,e_3^*)$, $(e_2,0)$ or $(0,e_3^*)$.
    \item If $(v,\alpha)$ is timelike, then it is in the $H$-orbit of $\sqrt{-\alpha (v)} (e_2,-e_2^*)$.
\end{itemize} 
\end{lemma}

In particular, $H$ acts transitively on unit speed space-like geodesics, and on unit speed time-like geodesics.

\subsubsection*{Cross ratios and distances}
We will now relate the geometry of $\mathbb{Y}$ to a cross ratio.
Let $(p_i,\ell_i)$ $i=1,2$ be a pair of points in $\mathbb{Y}$.

Suppose $p_i = [v_i]$ and $\ell_i = \ker \alpha_i$, for some vectors and functionals. 
We will define the cross ratio of the pair to be
$$CR((p_1,\ell_1),(p_2,\ell_2)) = \frac{\alpha_1 (v_2) \alpha_2 (v_1)}{\alpha_1 (v_1) \alpha_2 (v_2)}.$$
This number is independent of the choices of $v_i$ and $\alpha_i$, and is invariant under the $SL_3(\RR)$ action in pairs of points in $\mathbb{Y}$.

\begin{rmk}
      The value $CR((p_1, \ell_1), (p_2, \ell_2))$ coincides with the classical cross ratio associated to the points $p_1,p_2$, together with the intersections of the line $\overline{p_1 p_2}$ with $\ell_1$ and $\ell_2$, under the right normalization.
\end{rmk}

The following Lemma relates cross ratios with the pseudo-Riemannian geometry we introduced earlier:

\begin{lemma}
\label{lemma:crossratiogeometry}
    Let $(p,\ell) \in \mathbb{Y}$, and $X \in T_{(p,\ell)} \mathbb{Y}$ of length $\pm 1$ or zero, then
    $$CR((p,\ell), \exp (tX) \cdot (p,\ell)) = 
    \begin{cases}
        \cosh^2 (t),\: \text{if }X\text{ is spacelike},\\
        1,\: \text{if }X\text{ is lightlike},\\
        \cos^2 (t),\: \text{if }X \text{ is timelike}.
    \end{cases}$$
    Conversely, if $CR ((p_1,\ell_1), (p_2,\ell_2))> 0$, then the points $(p_1,\ell_1)$ and $(p_2,\ell_2)$ can be connected by a geodesic in $\mathbb{Y}$.
\end{lemma}

\begin{proof}
    The first item follows from an explicit computation with the vectors obtained in Lemma \ref{lemma:geodesicclassification}.
    The second item follows by noticing that, under a positive cross ratio assumption, we can send the point $p_1$, $p_2$, and $\ell_1 \cap \ell_2$ to the canonical basis and verify the claim there.
\end{proof}

As a consequence of the previous lemma, not every pair of parallel sets can be connected by a geodesic in $\mathbb{Y}$.
One can, for instance, take a look at pairs $(p_i,\ell_i)$ in non-degenerate relative position, meaning that $p_1$, $p_2$ and the intersections of $\overline{p_1 p_2}$ with $\ell_1$ and $\ell_2$ are distinct.

\begin{figure}[ht]
    \centering
    \includegraphics[width=0.8\textwidth]{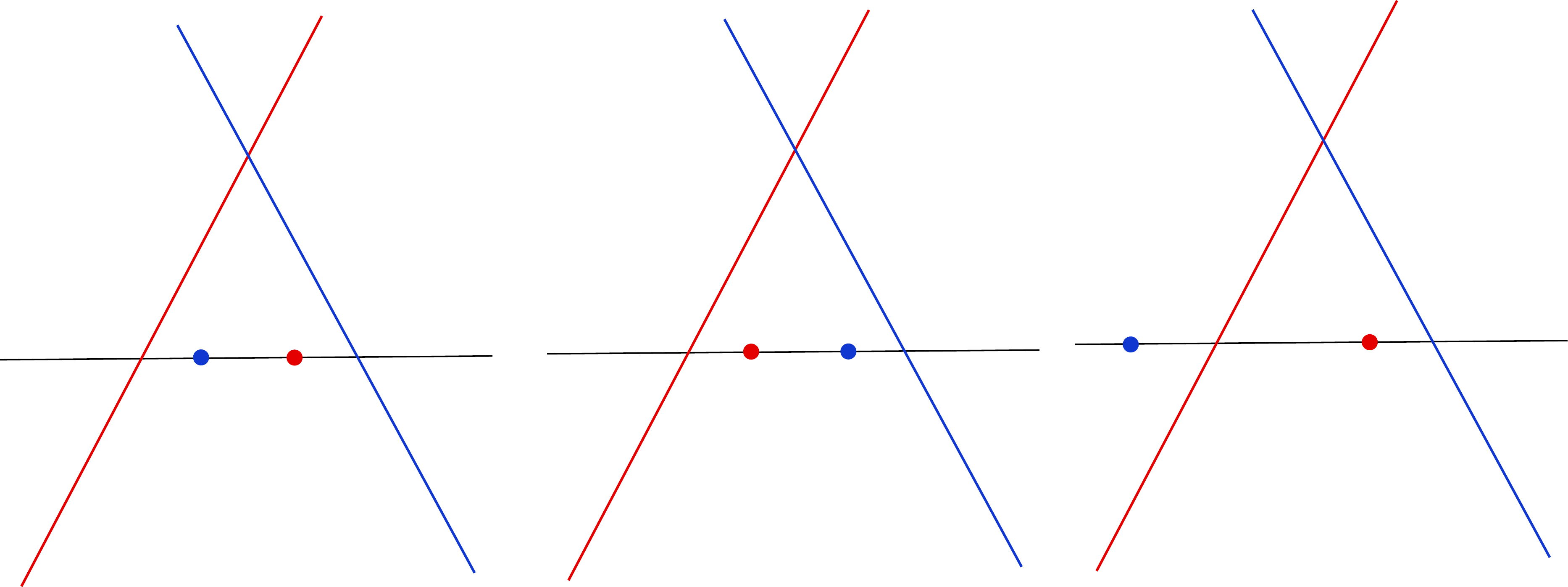}
    \caption{The relative positions of a pair in $\mathbb{Y}$. From left to right we call them (ND1), (ND2) and (ND3).}
    \label{domains:fig:relative-position}
\end{figure}

Figure \ref{domains:fig:relative-position} illustrates three such configurations, labeled (ND1), (ND2), and (ND3), respectively.
One can observe that (ND1) and (ND2) verify that $0 < CR < 1$ and $ CR> 1$ respectively.
On the other hand, (ND3) has a negative cross ratio, thus the pair of points and lines in this configuration cannot be joined by a geodesic.

One can also relate the cross ratio between a pair of points in $\mathbb{Y}$ to the trace of the composition of the associated involutions.

\begin{lemma}
\label{lemma:crossratioandtrace}
Let $\sigma_1$ and $\sigma_2$ be a pair of involutions associated to $(p_1,\ell_1),(p_2,\ell_2) \in \mathbb{Y}$, then
$$\operatorname{tr} (\sigma_1 \sigma_2) = 4 CR((p_1,\ell_1),(p_2,\ell_2)) -1.$$
\end{lemma}

\begin{proof}
Remark \ref{rmk:explicitinvolution} allows us to write $\sigma_i = \alpha_i \otimes v_i - Id$ where we normalize $\alpha_i (v_i) =2$ for $i=1,2$.
A direct computation shows that $\operatorname{tr} (\sigma_1 \sigma_2) = \alpha_1 (v_2) \alpha_2 (v_1) -1$.
\end{proof}

\begin{rmk}
Notice that $\sigma_1 \sigma_2$ always has $1$ as an eigenvalue, with eigenvector $\ell_1 \cap \ell_2$.
Thus the eigenvalues of $\sigma_1 \sigma_2$ are $\lambda$, $1$ and $\lambda^{-1}$, which can be computed using the formula above.
\end{rmk}

We finish the subsection by remarking that there is more to say about the relation of the cross ratio and the geometry of $\mathbb{Y}$.
More precisely, Labourie noticed that $\mathbb{Y}$ admits a natural symplectic structure that can be used to recover cross ratios (see \cite{labourie2007cross}).
From this, one can extract a more conceptual proof of Lemma \ref{lemma:crossratiogeometry}.

\subsection{Parallel sets}

We will now construct a correspondence between points in $\mathbb{Y}$ and totally geodesic \emph{parallel sets} on $\mathbb{X}$.

Suppose $\gamma$ is a singular geodesic in the symmetric space $\mathbb{X}$.
As explained in Subsection \ref{subsec:regular}, this implies that $\gamma' (0) \in T_{\gamma (0)} \mathbb{X}$ has a repeated eigenvalue, and thus has an associated pair of transverse eigenspaces $(p_\gamma,\ell_\gamma) \in \mathbb{Y}$ corresponding to the one and two-dimensional eigenspaces, respectively.
In this context, we say that $\gamma$ is \emph{asymptotic} to the pair $(p_\gamma,\ell_\gamma)$.

Unlike a regular geodesic, there is not a unique flat passing through $\gamma$.
A flat through $\gamma (0)$ is determined by a $\gamma(0)$-orthogonal basis.
It contains $\gamma$ if and only if the basis has the form $\{v_1,v_2,v_3\}$ with $v \in p_\gamma$ and $v_2,v_3 \in \ell_\gamma$
We define the \emph{parallel set} associated to the singular geodesic $\gamma$ to be
$$P_\gamma = \bigcup \{F : F \text{ flat, } \gamma \subset F\}.$$
Notice that there is a circle's worth of such flats, as two $q$-orthogonal bases with two vectors in $\ell_\gamma$ are related by an $SO(q|_{\ell\gamma}) \cong S^1$ transformation.

It is true in every symmetric space that parallel sets are totally geodesic submanifolds (see for instance \cite[Proposition 2.20.12]{eberlein1996geometry}).
We will explain this in our case with a self-contained argument that connects parallel sets to the space $\mathbb{Y}$.

As noted above, we can associate a singular geodesic with a point of $\mathbb{Y}$.
This point can, in turn, be associated to an involution $\sigma_{(p_\gamma,\ell_\gamma)}$, which is an isometry of the symmetric space.
It turns out that the fixed point set of this involution is precisely the parallel set $P_\gamma$.

\begin{prop}
\label{prop:singular=inv}
Let $\gamma$ be a singular geodesic of $\mathbb{X}$ asymptotic to a pair $(p_\gamma,\ell_\gamma) \in \mathbb{Y}$, and let $\sigma_\gamma$ be the involution associated with this pair.
Then the parallel set $P_\gamma$ coincides with the fixed point set of $\sigma_\gamma$, and it is a totally geodesic submanifold isometric to a product $\HH \times \RR$, where $\HH$ is a reducible slice.
Conversely, any totally geodesic copy of $\HH \times \RR$ is the parallel set of a singular geodesic.
\end{prop}

\begin{proof}
Notice that an involution $\sigma$ associated to $(p,\ell)$ fixes an inner product $q \in \mathbb{X}$ if and only if $\sigma \in SO(q)$.
This happens precisely when the orthogonal complement of $p$ coincides with $\ell$, thus
$$\operatorname{Fix} (\sigma) = { q \in \mathbb{X} : p^{\perp_q} = \ell}.$$
Clearly, the whole parallel set $P_\gamma$ is contained in $\operatorname{Fix} (\sigma_\gamma)$ by construction.
Conversely, any inner product fixed by $\sigma_\gamma$ restricted to $\ell_\gamma$ must share an orthogonal basis with $\gamma(0)|_\ell$.

Unimodular inner products that make the splitting $p_\gamma \oplus \ell_\gamma$ orthogonal are in bijection with (not necessarily unimodular) inner products on $\ell_\gamma$ by restriction.
The space of inner products on a two-dimensional vector space is isometric to $\mathbb{H} \times \mathbb{R}$, after identifying $\HH^2$ with the symmetric space of $PSL_2(\RR)$.

Finally, notice that any totally geodesic submanifold isometric to $\mathbb{H} \times \mathbb{R}$ verifies that ${z_0} \times \mathbb{R}$ is a geodesic contained in infinitely many flats of the form $\gamma \times \mathbb{R}$, for $\gamma$ a geodesic in $\mathbb{H}$ passing through $z$.
Thus, this geodesic is singular.
\end{proof}

Therefore, the space $\mathbb{Y}$ is in bijection with the set of \emph{parallel sets} of $\mathbb{X}$.
We will now see that the relative position of pairs of parallel sets reflects on the relative position of the corresponding points of $\mathbb{Y}$ in terms of its geometry.

\begin{prop}
\label{prop:relativeposition}
Let $P_1,P_2$ be the parallel sets associated to $(p_1,\ell_1),(p_2,\ell_2)\in\mathbb{Y}$, let $\sigma_1,\sigma_2$ be the corresponding involutions, and suppose that $CR((p_1,\ell_1),(p_2,\ell_2)) \notin \{0,1\}$.
Then $d_{\mathbb{X}}(P_1,P_2)$ is realized by a pair of points in $\mathbb{X}$, and 
\begin{enumerate}
\item $P_1 \cap P_2 \neq \emptyset$ if and only if $0 < CR < 1$. In this case, $(p_1,\ell_1)$ and $(p_2,\ell_2)$ are joined by a timelike geodesic.
\item If $CR > 1$, then $(p_1,\ell_1)$ and $(p_2,\ell_2)$ are joined by a spacelike geodesic, and $CR = \cosh^2 d_\mathbb{X} (P_1,P_2)$.
\item If $CR < 0$, then $(p_1,\ell_1)$ and $(p_2,\ell_2)$ are not joined by a geodesic, and $CR = -\sinh^2 d_\mathbb{X} (P_1,P_2)$.
\end{enumerate}
\end{prop}

\begin{proof}
By Lemma \ref{lemma:crossratioandtrace}, the eigenvalues of $\sigma_1 \sigma_2$ are $1,\lambda,\lambda^{-1}$ with $\lambda + \lambda^{-1} = 4CR - 2$.
Since $CR \notin \{0,1\}$, these are distinct, so $\sigma_1 \sigma_2$ is diagonalizable over $\mathbb{C}$.

If $\lambda$ is complex, it is necessarily in the unit circle, and thus $\sigma_1 \sigma_2$ preserves some inner product $q \in \mathbb{X}$.
Otherwise, $\lambda$ is real, $\sigma_1 \sigma_2$ is diagonal over $\mathbb{R}$, and preserves a flat $F_{\sigma_1 \sigma_2}$ in the symmetric space.
In any case, the translation distance of $\sigma_1 \sigma_2$ is $\log |\lambda|$ (with our normalization of the metric on $\mathbb{X}$).
Given $C$ the set of points with minimal displacement:
$$C = \{ q \in \mathbb{X}: d(q,\sigma_1 \sigma_2 q) = \log |\lambda|\},$$
we saw that $C$ is non-empty in both cases, and it is convex, as it is either a flat, or the set of inner products making $\sigma_1 \sigma_2$ orthogonal.

We claim that $d_\mathbb{X} (P_1,P_2) = \log |\lambda|/2$, and that it is attained.
For $x_1\in P_1$ and $x_2\in P_2$ we have
$$d(x_2, \sigma_1 \sigma_2 x_2) = d(x_2,\sigma_1 x_2) \leq d(x_2,x_1) + d(x_1,\sigma_1 x_2) = 2 d(x_1,x_2),$$
so $\log |\lambda| \leq 2d_\mathbb{X} (P_1,P_2)$.
Conversely, $\sigma_1$ conjugates $\sigma_1 \sigma_2$ with its inverse, so $\sigma_1$ preserves $C$.
The midpoint of $y$ and $\sigma_1 y$, for any $y \in C$, is then a point $x_1 \in P_1 \cap C$.
Since $d(x_1,\sigma_2 x_1) = d(\sigma_1 x_1, \sigma_1\sigma_2 x_1) = \log |\lambda|$, the midpoint $m$ of $[x_1,\sigma_2 x_1]$ is fixed by $\sigma_2$, so $m\in P_2$ and $d_\mathbb{X} (P_1,P_2) \leq d(x_1,m) = \log |\lambda|/2$.

Therefore $\log|\lambda| = 2d_\mathbb{X} (P_1,P_2)$.
If $0<CR<1$, then $\lambda$ lies on the unit circle, and $P_1\cap P_2 \neq\emptyset$.
If $\lambda>0$, then $4CR - 2 = 2\cosh 2 d_{\mathbb{X}} (P_1,P_2)$, that is, $CR = \cosh^2 d_\mathbb{X} (P_1,P_2)$.
The case where $\lambda < 0$ is identical.
The statements about geodesics follow from Lemma \ref{lemma:crossratiogeometry}.
\end{proof}

The cases where $CR \in \{0,1\}$ are somewhat degenerate, and we leave it for the reader to analyze.
We note that the $CR = 1$ scenario includes the case of lightlike geodesics, which correspond to asymptotical parallel sets.

It is also true that points in $\mathbb{X}$ define totally geodesic submanifolds of $\mathbb{Y}$, though we will not explore this here.

\begin{lemma}
\label{lemma:convergencia-paralelos}
Let $(p_n, \ell_n)$ be a sequence in $\mathbb{Y}$ associated to parallel sets $(P_n)$, and assume that there exists a convergence sequence $(x_n) \in \mathbb{X}^\mathbb{N}$ such that $x_n \in P_n$. 
Then $(p_n, \ell_n)$ has a convergent subsequence in $\mathbb{Y}$.
\end{lemma} 

\begin{proof}
Observe that each $(x_n)$ has a distinguished singular direction given by some $v_i \in T_{x_i} \mathbb{X}$ so that $P_i$ is the parallel set associated to the singular geodesic in that direction.
Since the sequence $(x_n)$ converges, we can take a convergent subsequence for $(v_i)$ in $T \mathbb{X}$. 
Moreover, since the set of singular directions is closed, any limit of direction $v \in T\mathbb{X}$ defines a parallel set $P$ with the desired properties.
\end{proof}

\subsection{Irreducible slices and orthogonal parallel sets}
\label{subsec:irredparallel}

In Subsection \ref{subsub:hypslices}, we studied hyperbolic slices, which are defined from bilinear forms. We will study the relationship between these and the parallel sets parametrized by $\mathbb{Y}$.
The key remark is that the closest point projection $\pi: \mathbb{X} \rightarrow \mathcal{H}$ onto a slice $\mathcal{H}$ defines a fibration of $\mathbb{X}$ with parallel sets as fibers as hinted by Figure \ref{fig:walls}.
We make this formal in the following proposition, which motivates most of this paper:

\begin{prop}
\label{prop:slicefibration}
Let $B$ be a symmetric bilinear form of signature $(1,2)$ with unit determinant defining an ellipse $\mathcal{E}_B = \mathbb{P}(\{v : B(v,v) > 0\})$, and $\pi : \mathbb{X} \to \mathcal{H}_B$ be the closest point projection.
For $p \in \mathcal{E}_B$, let $\sigma_p$ be the involution associated to $(p,p^{\perp_B}) \in \mathbb{Y}$, and $P_p = \operatorname{Fix}(\sigma_p)$ its parallel set.
Then
\begin{enumerate}
    \item the map $p \mapsto q_p = B(\sigma_p \cdot,\cdot)$ is a bijection $\mathcal{E}_B \to \mathcal{H}_B$.
    \item $T_{q_p} P_p = (T_{q_p}\mathcal{H}_B)^\perp$ and $\pi^{-1}(q_p) = P_p$. 
    In particular, $\pi$ is a fibration whose fibers are parallel sets.
    \item The map $\Phi : \mathcal{H}_B \to \mathbb{Y}$, $q_p \mapsto (p,p^{\perp_B})$, is an isometric embedding onto a totally geodesic spacelike surface.
\end{enumerate}
\end{prop}

\begin{proof}
Let $q \in \mathbb{X}$, recall that we denote by $[B]_q$ the matrix such that $B(\cdot, \cdot) = q([B]_q\cdot, \cdot)$.
Note that $[B]_q$ is $q$-self-adjoint, so $q \in \mathcal{H}_B$ (i.e. $[B]_q \in O(q)$), if and only if $[B]_q^2 = Id$, that is, if it is a $q$-orthogonal involution.

\begin{enumerate}
\item Let $q \in \mathcal{H}_B$ and $(p,\ell)$ be the eigenspaces of $[B]_q$.
For $v \in p$ and $w \in \ell$ we have
$$B(v,v) = q(v,v) > 0, \quad B(w,w) = -q(w,w) < 0,$$ $$\quad B(v,w) = q(v,w) = 0.$$
Since $B$ has signature $(1,2)$, this forces $\dim p = 1$, $p \in \mathcal{E}_B$ and $\ell = p^{\perp_B}$.
Hence $[B]_q = \sigma_p$ and $q = B(\sigma_p \cdot,\cdot) = q_p$.
Conversely, the same three identities show that $q_p$ is positive definite.
The map is injective because $p$ is precisely the line fixed by $\sigma_p$.

\item Again let $q \in \mathcal{H}_B$. By Remark \ref{rmk:compinvolutions} and the item above, we have 
$I_B I_{q_p} = \sigma_p$. We have the following formulas:
$$d_{q_p}I_B|_{T_{q_p} \mathcal{H}_B} = Id, \qquad d_{q_p} I_{q_p} = -Id$$
Moreover, because $\mathcal{H}_B$ is the fixed point set of $I_B$, we have that $d_{q_p}I_B$ is $-Id$ on $(T_{q_p} \mathcal{H}_B)^\perp$.

So $d_{q_p} \sigma_p$ is $-Id$ on $T_{q_p}\mathcal{H}_B$ and $Id$ on $(T_{q_p} \mathcal{H}_B)^\perp$, and therefore $T_{q_p} P_p = (T_q\mathcal{H}_B)^\perp$.
Notice that $\pi(x) = q$ if and only if the geodesic from $q$ to $x$ is orthogonal to $\mathcal{H}_B$ at $q$.
Hence $\pi^{-1}(q) = P_p$.

\item The map $\Phi$ is $SO(B)$-equivariant, so for $X \in T_q\mathcal{H}_B$ it sends the geodesic $e^{tX}\cdot q$ to $e^{tX}\cdot(p,p^{\perp_B})$.
Since $X$ is in the $-1$ eigenspace of $\operatorname{Ad} (\sigma_p)$, the latter is a geodesic of $\mathbb{Y}$ with initial velocity $X$.
Both metrics are $\frac{1}{2}\operatorname{tr}(X^2)$, which is positive on symmetric matrices, so $\Phi$ is a spacelike isometric immersion.
Injectivity follows from the first item.
\end{enumerate}
\end{proof}

In particular, the preimage of a point under the closest point projection $\pi: \mathbb{X} \rightarrow \mathcal{H}$ is a parallel set.

There is a more geometric way of thinking of this construction.
We can define a hyperboloid $\mathcal{C}_B = \{v \in \mathbb{R}^3: B(v,v) = 1\}$.
The tangent space to $\mathcal{C}_B$ at a point $v$ is precisely the $B$-orthogonal complement of the line defined by $v$, which is a negative definite subspace.
In particular, $-B|_{T \mathcal{C}_B}$ is an inner product.
We can extend it to an inner product in $\mathbb{R}^3$ by extending it as $B|_{\langle v \rangle}$.
This extension associated to $B$ is exactly $q_{\langle v \rangle}$.

As in the previous section, one can prove the proposition more generally.
More precisely, given a reflective totally geodesic submanifold in a symmetric space, the orthogonal of the tangent space at a point is again a reflective submanifold.

\subsection{Spacelike geodesics revisited}
\label{subsection:spacelike}

Proposition \ref{prop:slicefibration} gives us a synthetic description of spacelike geodesics.
More precisely, pulling back geodesics in $\mathcal{H}_B$ to $\mathcal{E}_B$ gives a segment in $\mathcal{E}_B$, corresponding to a geodesic in the Klein model of hyperbolic space.
The trace of the corresponding geodesic in $\mathbb{Y}$ corresponds to taking points on this line and their orthogonal complement.
Any spacelike geodesic has this form, since two spacelike geodesics lie in the same $G$-orbit.

We can compactify the space $\mathbb{Y}$ by considering it inside $\mathbb{RP}^2 \times (\mathbb{RP}^2)^*$, and think of the space of full flags $\mathcal{F}$ as its boundary.
Given $\gamma : \mathbb{R} \rightarrow \mathbb{Y}$, a geodesic obtained from an ellipse $B$, we see that $\lim_{t \to \pm \infty} \gamma (t) \in \mathcal{F}$ is a pair of flags tangent to the ellipse $\partial \mathcal{E}_B$ at the endpoints of the segment defined by the geodesic.
Moreover, these two flags are in general position.

\begin{figure}
    \centering
    \includegraphics[scale=.2]{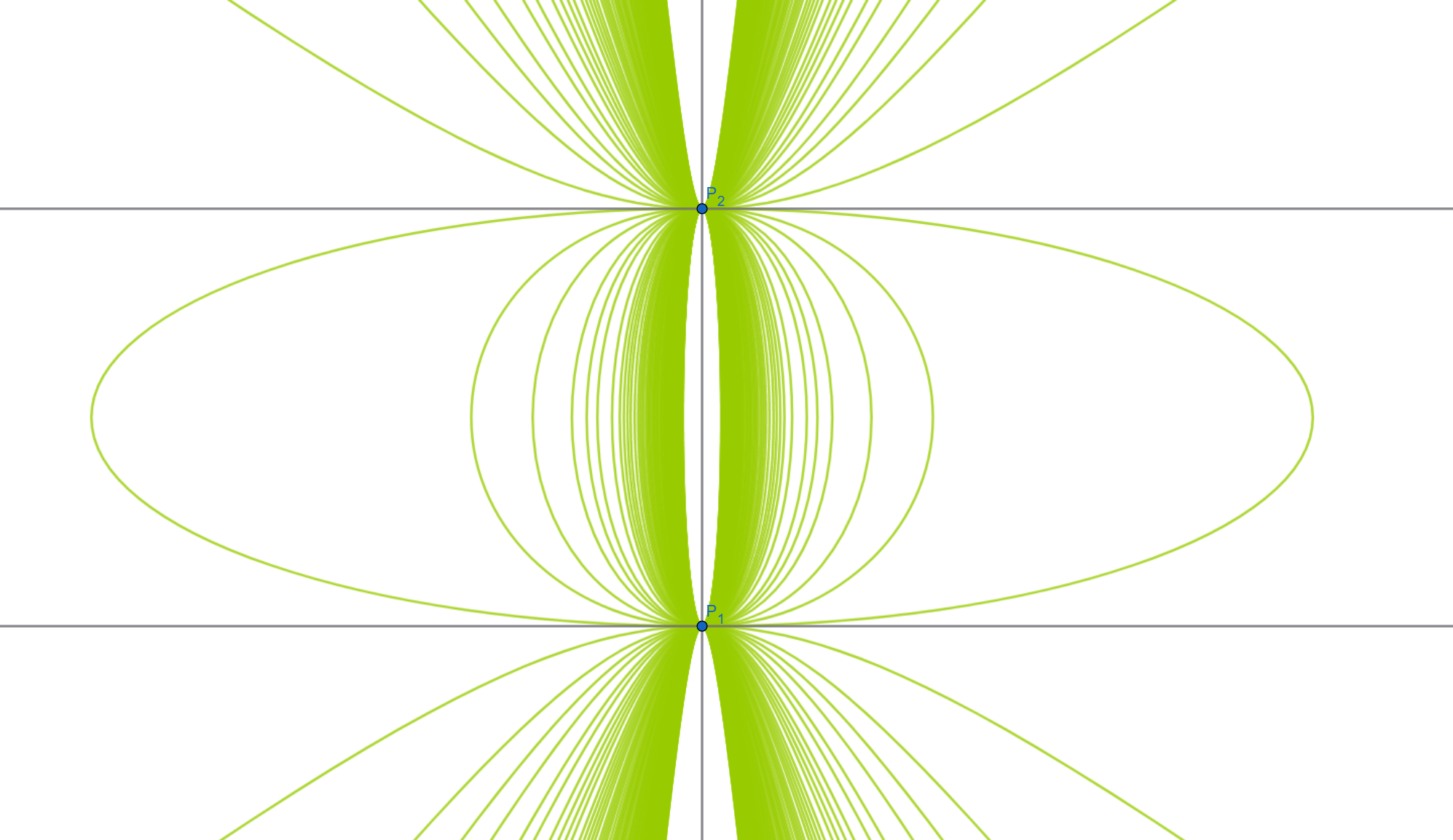}
    \caption{Pencil of conics tangent to $F_1$ and $F_2$}
    \label{fig:pencil}
\end{figure}

Conversely, we can explicitly compute the geodesics asymptotic to a pair $(F_1,F_2)$. 
First, we observe that two ellipses tangent to a given pair of flags, and sharing a segment joining the points of the flag define the same spacelike geodesics.

\begin{lemma}
\label{lemma:polar}
Let $F_i = (p_i,\ell_i)$, $i=1,2$, be transverse flags, and choose $v_i \in p_i$ and $\alpha_i$ with $\ker\alpha_i = \ell_i$ such that $\alpha_1(v_2) = \alpha_2(v_1) = 1$.
Let $B$ be a form of signature $(1,2)$ whose conic is tangent to $F_1$ and $F_2$, and let $v \in p_1\oplus p_2$ with $B(v,v) \neq 0$.
Then $[v]^{\perp_B} = [\eta_v]$, where
$$\eta_v = \alpha_1(v)\,\alpha_2 + \alpha_2(v)\,\alpha_1.$$
In particular, $[v]^{\perp_B}$ does not depend on $B$.
\end{lemma}

\begin{proof}
Since $G$ acts transitively on pairs of transverse flags and the formula for $\eta_v$ is equivariant, we may assume $F_1 = ([e_1],[e_3^*])$, $F_2 = ([e_3],[e_1^*])$.
Tangency means $B(e_1,\cdot)$ and $B(e_3,\cdot)$ are proportional to $e_3^*$ and $e_1^*$ respectively, so up to scaling $B$ is an element of the family $B_\lambda(x,y,z) = \lambda xz - y^2$ for some $\lambda \neq 0$ (see Figure \ref{fig:pencil}).
For $v = e_3 + s e_1$ we get $B_\lambda(v,\cdot) = \frac{\lambda}{2}(e_1^* + s e_3^*) = \frac{\lambda}{2}\,\eta_v$.
\end{proof}

The following gives us an explicit expression for a spacelike geodesic asymptotic to a pair of flags in general position.

\begin{prop}
\label{prop:asymptoticgeodesics}
Let $F_1, F_2$ and $v_i,\alpha_i$ be as in Lemma \ref{lemma:polar}.
Let $c$ and $c'$ be the connected components of $\mathbb{RP}^2 \setminus (\ell_1 \cup \ell_2)$, with $[v_2 + s v_1] \in c$ for $s>0$. 
Then
\begin{enumerate}
\item the curves
$$\gamma_c (t) = ( [v_2 + e^{2t}v_1], [\alpha_2 + e^{2t}\alpha_1]), \: \gamma_{c'}(t) = ([v_2 - e^{2t}v_1], [\alpha_2 - e^{2t}\alpha_1]),$$
are unit speed spacelike geodesics of $\mathbb{Y}$, converging to $F_2$ as $t\to-\infty$ and to $F_1$ as $t\to+\infty$.
They are the only spacelike geodesics of $\mathbb{Y}$ with these limits.
\item For every $t,t'\in\RR$,
$$CR (\gamma_c(t),\gamma_c(t')) = \cosh^2(t-t'), \: CR (\gamma_c(t),\gamma_{c'}(t')) = -\sinh^2(t-t').$$
In particular, the parallel sets of $\gamma_c(t)$ and $\gamma_{c'}(t')$ are at distance $|t-t'|$, and they intersect when $t = t'$.
\end{enumerate}
\end{prop}

\begin{proof}
Writing $v_s = v_2 + s v_1$, Lemma \ref{lemma:polar} gives $\eta_{v_s}(v_{s'}) = s + s'$, and $\eta_{v_s}(v_s) = 2s$, thus:
$$CR ((v_s,\eta_{v_s}),(v_{s'},\eta_{v_{s'}})) = \frac{(s+s')^2}{4ss'}.$$
Substituting $s = \pm e^{2t}$ and $s' = \pm e^{2t'}$ gives the formulas in (2).
The limits in (1) are immediate.

Choose $B$ in the pencil of Lemma \ref{lemma:polar} whose ellipsoid $\mathcal{E}_B$ contains the segment $s= \{[v_2 + s v_1]: s>0\}\subset c$.
By Proposition \ref{prop:slicefibration}, $\{(p,p^{\perp_B}) : p \in s\}$ is the image of a geodesic of $\mathcal{H}_B$, hence a spacelike geodesic of $\mathbb{Y}$.
The same argument applies for $\gamma_{c'}$.
These geodesics are different because when $t = t'$, $\gamma_c (t)$ intersects $\gamma_c (t')$, but are different parallel set.
On the other hand, Lemma \ref{lemma:polar} tells us that any spacelike geodesic asymptotic to $F_1$ and $F_2$ coincides with $\gamma_c$ or $\gamma_{c'}$.
\end{proof}

\subsection{The Walls} 
\label{section:medialwall}

We will now use the duality between $\mathbb{X}$ and $\mathbb{Y}$ to define a special type of hypersurface in $\mathbb{X}$ that we call a \emph{wall}.
Each parallel set is three-dimensional; thus, a one-parameter family of disjoint parallel sets is automatically codimension one.

As seen in the previous section, any spacelike geodesic in $\mathbb{Y}$ is asymptotic to a pair of transverse flags $(F_1,F_2)$.
Conversely, any pair of transverse flags is asymptotic to precisely two spacelike geodesics.
A unique geodesic is specified once we fix a favorite connected component $c$ of $\RR \PP^2 \setminus (l_1 \cup l_2)$ containing the points defined by the geodesic in $\mathbb{Y}$.

\begin{defi}
\label{defi:walls}
Given a pair of transverse flags $(F_1,F_2)$ and a conected component $c$ as above, we define its associated \emph{wall} as:
$$\mathcal{W}(F_1,F_2,c) = \bigcup_{t \in \mathbb{R}} P(t),$$
where $P(t)$ are parallel sets defined by the spacelike geodesic joining $F_1$ and $F_2$ corresponding to the component $c$.
\end{defi}

We sometimes denote $\mathcal{W}(F_1, F_2)$ when $c$ is uniquely determined.
By everything we did before, we observe that parallel sets in the definition of $\mathcal{W}(F_1,F_2,c)$ are pairwise disjoint.
Since $\mathbb{X}$ is five-dimensional, the walls are hypersurfaces.

The following fact will be useful later on.

\begin{rmk}
\label{rmk:finitelymanyint}
Any wall $\mathcal{W} \subset \mathbb{X}$ is an analytic hypersurface.
In particular, a compact geodesic segment in $\mathbb{X}$ is either contained in $\mathcal{W}$ or intersects it finitely many times.
\end{rmk}

Any wall can be associated to a one-parameter family of irreducible hyperbolic slices, corresponding to ellipses tangent to $F_1$ and $F_2$ and contained in $c$. 
Pick any $\mathcal{H}$ in this family, and consider $\pi : \mathbb{X} \rightarrow \mathcal{H}$ the closest point projection, then $\mathcal{W}(F_1,F_2,c)$ is the preimage of a geodesic $\gamma$ asymptotic to $F_1$ and $F_2$.
Particularly, we get the following.

\begin{lemma}
\label{lemma:walls-separate}
Walls $\mathcal{W}(F_1,F_2,c)$ separate the symmetric space into two connected components, each isometric to $\pi^{-1} (H)$, for $\pi: \mathbb{X} \rightarrow \mathcal{H}$ the projection onto an irreducible hyperbolic slice, and $H \subset \mathcal{H}$ a half-space.
\end{lemma}

Inspired by \cite[Chapter 4]{goldman1999complex}, we will call the parallel sets $P(t)$ the \textit{slices} of the wall $\mathcal{W}(\gamma,\mathcal{H})$.

The wall construction is natural, in the following sense:

\begin{lemma}
\label{lemma:naturality}
Suppose $g \in G$ is an element sending a triple $(F_1,F_2,c)$ as before to $(F_1',F_2',c')$, then $g\mathcal{W}(F_1,F_2,c) = \mathcal{W}(F_1',F_2',c')$.
\end{lemma}

\begin{proof}
Any such $g$ sends an ellipsoid tangent to $F_1$ and $F_2$ and contained in $c$ to an ellipsoid tangent to $F_1'$ and $F_2'$ and contained in $c'$.
\end{proof}

%% file: sections/triples.tex
\section{Parallel domains}
\label{sec:parallel}

This section constructs a domain $\mathcal{D}(F_1, F_2, F_3)$ on the symmetric space $\mathbb{X}$ associated with a positive triple of flags $(F_1, F_2, F_3)$ with boundary components consisting of walls.
These domains mimic how we associate to a triple of points in $\partial \mathbb{H}^2$ a unique ideal triangle in hyperbolic space.

\subsection{Triples of flags in general position}
\label{domains:subsec:projinv}

In this section, we briefly review the notion of a triple of flags in general position, and describe a projective invariant called the \emph{triple ratio} that parametrizes the $SL_3 (\mathbb{R})$-orbits of triples of flags in general position.

Whenever $F \in \mathcal{F}$ is a flag, we will denote by $p$ and $\ell$ the corresponding point and line.
Recall that $\mathcal{F}^{(2)}$ denotes the set of pairs of flags in general position, as introduced in Section \ref{subsec:regular}.

\begin{defi}
We say that a triple of flags $(F_i)_{i=1,2,3} \in \mathcal{F}^{3}$ is in \emph{general position} if every pair is in general position, the three corresponding points are not aligned, and the three lines are not concurrent at a point.
We denote the set of triples in general position as $\mathcal{F}^{(3)}$.
\end{defi}

As mentioned before, $SL_3 (\mathbb{R})$-orbits in $\mathcal{F}^{(3)}$ is parametrized by a unique real number.
Given a triple $(F_i) \in \mathcal{F}^{(3)}$ in general position, we define their \emph{triple ratio} as:
$$TR(F_1, F_2, F_3) = \frac{\alpha_1 (v_2) \alpha_2 (v_3) \alpha_3 (v_1)}{\alpha_1 (v_3) \alpha_2 (v_1) \alpha_3 (v_2)},$$
where $(v_i,\alpha_i) \in \RR^3 \times (\RR^3)^*$ are vectors such that $[v_i] = p_i$, $\ker \alpha_i = \ell_i$.
The definition is independent of any choice.
This quantity was originally introduced in \cite[Section 2]{fock2007moduli}.

We say that a triple of flags is \textit{positive} (resp. negative) if they have a positive (resp. negative) triple ratio. 
Figure \ref{fig:triplet-tr} in the introduction depicts both a positive and a negative triple.
A triple is positive if and only if there is an affine chart in which the convex hull of the flag points lies in a triangle bounded by the three lines.

The triple ratio is a complete $SL_3 (\mathbb{R})$-invariant for triples in general position in the following sense.

\begin{lemma}
\label{lemma:canonicalpos}
Two triples $(F_i),(F_i') \in \mathcal{F}^{(3)}$ are in the same $SL_3 (\mathbb{R})$-orbit if and only if $TR(F_1,F_2,F_3) = TR (F_1',F_2',F_3')$.
\end{lemma}

\begin{proof}
Notice that since the lines of a triple $(F_i) \in \mathcal{F}^{(3)}$ are pairwise different and not concurrent at a point, we can assume that $\ell_i = \ker e_i^*$, for $\{e_i\}_{i=1,2,3}$ the canonical basis of $SL_3 (\mathbb{R})$.
Up to further composing with a positive diagonal matrix, we can assume that
$$F_1 (s) = ([s e_2 + e_3], [e_1^*]),\: F_2 (s) = ([e_1 + s e_3], [e_2^*]), \: F_3 (s) = ([se_1 + e_2], [e_3^*]),$$
for some $s \in \mathbb{R}$.

Note that $s \neq 0$ and $-1$, since the points cannot lie at the intersection of the lines or be collinear.
A direct computation shows that $TR(F_1(s), F_2(s), F_3(s)) = s^{-3}$, thus the triple ratio determines the $s$ uniquely.
\end{proof}

As a consequence of the proof, triple ratios of a triple in $\mathcal{F}^{(3)}$ cannot attain the values $0$ and $-1$.

\begin{rmk}
Notice that the expression for the triple ratio is well defined as long as pairs of flags are in general position.
For our purposes, we wanted $\mathcal{F}^{(3)}$ to avoid the scenario in which the lines are concurrent at a point, or the three points are aligned.
The latter configuration corresponds to a triple ratio of $-1$, as Menelaus' theorem predicts.
\end{rmk}

There are two natural extra symmetries of a triple of flags.
One of them is a \emph{duality} in the sense that we now introduce.
This has been observed by Schwartz in \cite[Section 5]{schwartz2025leretourpappus}.

Suppose $B$ is a nondegenerate, symmetric bilinear form in $\mathbb{R}^3$.
Such a form defines a bijection $D_B$ of $\mathbb{RP}^2 \times (\mathbb{RP}^2)^*$ into itself by declaring
$$D_B (p,\ell) = (\ell^{\perp_B}, p^{\perp_B}),$$
which we call the \emph{duality associated to $B$}.
A duality is a transformation coming from this procedure for some $B$.
Observe that $D_B^2 = Id$ by construction.

The following lemma collects elementary computations that will be useful later.

\begin{lemma}
\label{rmk:composition-dualities}
Let $B$ be a non-degenerate, symmetric bilinear form of unit determinant and with signature $(1,2)$.
\begin{enumerate}
\item Given $B'$ another such form, then the composition $D_{B'} \circ D_{B}$ coincides with multiplication by $[B']^{-1}\in SL_3 (\mathbb{R})$ defined via $B' (v,w) = B ([B']v,w)$.
\item If $P \subset \mathbb{X}$ is a parallel set associated to a pair $(p,\ell) \in \mathbb{Y}$, then $I_{B} (P)$ is a parallel set associated to $D_{B_1} (p,\ell)$.
\item The duality $D_B : \mathbb{Y} \rightarrow \mathbb{Y}$ is an isometry for the pseudo-Riemannian metric on $\mathbb{Y}$.
\end{enumerate}
\end{lemma}

\begin{proof}
\phantom{v}

\begin{enumerate}
    \item Follows by direct computation. 
    \item Follows from the fact that if a parallel set is defined by a tangent direction $X \in \mathfrak{p}_q$ with eigenspaces $p$ and $\ell$, then $X^{*_B}$ has eigenspaces $\ell^{\perp_B}$ and $p^{\perp_B}$.
    \item\label{it:composition-dualities3} Follows from the fact that $I_B$ is an isometry, and distances between parallel sets in the symmetric space are related to distances in $\mathbb{Y}$.
\end{enumerate}
\end{proof}

Dualities complete the list of symmetries of a triple of flags.

\begin{lemma}
\label{lemma:duality}
Let $(F_i) \in \mathcal{F}^{(3)}$ be a triple in general position.
Taking indices mod $3$, we get:
\begin{enumerate}
\item There exists a unique $S \in SL_3 (\mathbb{R})$ such that $S F_i = F_{i+1}$.
\item There exists a unique duality $D_i$ such that $D_i F_i = F_i$, and $D_i F_{i+1} = F_{i-1}$.
\item $D_i \circ D_{i+1} = S$ for any $i=1,2,3$.
\end{enumerate}
\end{lemma}

\begin{proof}

$(1)$ This follows from the fact that $TR(F_1,F_2,F_3) = TR(F_2,F_3,F_1)$ and Lemma \ref{lemma:canonicalpos}.
This transformation is unique: if $S$ and $S'$ are two such transformations, then $S' S^{-1}$ fixes four points in general position and therefore acts as the identity. \\
$(2)$ Assume that the triple $(F_i)$ is of the form $F_i (s)$ as in Lemma \ref{lemma:canonicalpos}.
A bilinear form $B$ describing a duality $D_3$ as in the statement should satisfy that:
\begin{itemize}
    \item $s e_2 + e_3$ is $B$-orthogonal to $[e_1] \oplus [e_3]$.
    \item $e_1 + s e_3$ is $B$-orthogonal to $[e_2] \oplus [e_3]$.
    \item $e_3$ is $B-$orthogonal to $[s e_2 + e_3] \oplus [e_1 + se_3 ]$.
    \item The $B-$orthogonal complement of $s e_1 + e_2$ is precisely $[e_1] \oplus [e_2]$.
\end{itemize}
Putting those conditions together, we obtain that defining
$$[B] =
\begin{pmatrix}
1 & -s & s^2\\
-s & s^2 & 1\\
s^2 & 1 & -s
\end{pmatrix},$$
then $B(v,w) = q_0 ([B]v,w)$ verifies the desired properties, for $q_0$ the standard inner product.
Uniqueness follows from Lemma \ref{rmk:composition-dualities} and the argument above.\\
Item \ref{it:composition-dualities3} of Lemma \ref{rmk:composition-dualities} shows that $D_i \circ D_{i+1}$ cyclically permutes $(F_i)$, and is an element of $PGL_3 (\mathbb{R})$.
By the first point, this element coincides with $S$.
\end{proof}

\subsection{Definition of the domains}
In this section, we introduce the notion of a \emph{parallel domain} associated with a positive triple of flags, as promised in the introduction.

Recall that a triple of flags in general position $(F_1, F_2, F_3) \in \mathcal{F}^{(3)}$ is positive if and only if the points corresponding to the flags bound a triangle whose interior $T$ is entirely contained in a connected component of the complement of the lines of $F_i$ (see Figure \ref{fig:triplet-tr}).
This lets us choose, for each $i$, a canonical connected component $c_i$ of the complement of the lines of $F_{i-1}, F_{i+1}$ containing $T$.
Therefore, for a positive triple, we can assign three walls $\mathcal{W}_i = \mathcal{W} (F_{i-1}, F_{i+1},c_i)$.

Crucially, these three walls are disjoint.

\begin{prop}
\label{prop:criterion_disjointness}
Let $(F_i) \in \mathcal{F}^3$ be a positive triple of flags, and $c_i$ the connected component of $\RR\PP^2 \setminus \ell_{i+1} \cup \ell_{i - 1}$ containing $p_i$, and $\mathcal{W}_i = \mathcal{W}(F_{i + 1}, F_{i - 1}, c_i)$ as above. 
Then $\mathcal{W}_i \cap \mathcal{W}_{i+1} = \emptyset$ for any $i = 1,2,3$.
\end{prop}

\begin{figure}[ht]
    \centering
    \includegraphics[width=0.8\textwidth]{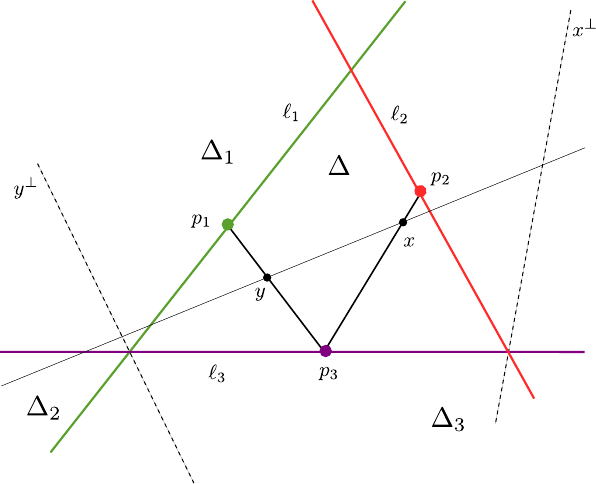}
    \caption{A representation of Proposition \ref{prop:criterion_disjointness}. The pairs $(x, x^\perp)$ and $(y, y^\perp)$ are in spacelike position, and therefore the parallel sets defining them are disjoint.}
    \label{fig:walls-disjuntas}
\end{figure}

\begin{proof}
 See Figure \ref{fig:walls-disjuntas} for a representation of this argument. 
 Let $\Delta$ be the projective triangle containing the three points of the flag triple on its boundary, and let $\Delta_i$ be the triangle determined by the three lines that contains only the point $p_i$ on its boundary. 
    
 Without loss of generality, we will prove that $\mathcal{W}_1$ and $\mathcal{W}_2$ are disjoint. 
 Let $\gamma_1$ be the segment in $c_1$ joining $p_2$ with $p_3$, and $\gamma_2$ the one joining $p_1$ and $p_3$, and pick $x \in \gamma_1$ and $y\in \gamma_2$.
There are parallel sets in $\mathcal{W}_1$ and $\mathcal{W}_2$ whose points coincide with $x$ and $y$ respectively.
We denote by $x^\perp$ and $y^\perp$ the lines of the corresponding parallel sets.

Let $\gamma_{x, y}$ be the line joining $x$ and $y$, oriented from $x$ to $y$. 
This line will follow one of three possible itineraries: 
$$\Delta \rightarrow \Delta_1 \rightarrow \Delta_2 \rightarrow \Delta, \qquad \Delta \rightarrow \Delta_3 \rightarrow \Delta_2 \rightarrow \Delta, \qquad \Delta \rightarrow  \Delta_1 \rightarrow \Delta_3 \rightarrow \Delta.$$
The proof is analogous for any of the three possible itineraries, so we will prove it for the first.
Notice that $\gamma_{x, y}$ can only go from $\Delta$ to $\Delta_1$ through $\ell_1$, and from $\Delta_1$ to $\Delta_2$ from $\ell_3$. 
 Note that $y^\perp$ is a line passing through $\ell_1 \cap \ell_3$ contained in $\Delta_1 \cup \Delta_3$, and the segment of $\gamma_{x, y}$ included in $\Delta_1$ intersects all such lines. 
 This implies that the segment cuts $y^\perp$ in $\Delta_1$. Together with an analogous argument for $x^\perp$, this implies that the pair $(x, x^\perp), (y, y^\perp)$ in $\mathbb{Y}$ are in position (ND2) (i.e., spacelike), as seen in Figure \ref{domains:fig:relative-position}, which in turn implies that the walls are disjoint.
\end{proof}

It is therefore natural to assume that there is a common region bounded by these three walls.
We prove it in the following lemma.

\begin{lemma}
\label{lemma:pingpongregions}
Given $\mathcal{W}_i$, $i=1,2,3$ associated with a positive triple $(F_1,F_2,F_3)$ as above, then $\mathbb{X} \setminus \bigcup_{i=1}^3 \mathcal{W}_i$ is a union of four connected components $\mathcal{D} \cup \Delta_1 \cup \Delta_2 \cup \Delta_3$, where
\begin{itemize}
\item $\partial \mathcal{D} = \bigcup_{i=1}^3 \mathcal{W}_i$ and $\mathcal{D}$ is invariant under the order three involution $S$ that cyclically permutes $F_i$.
\item $\partial \Delta_i = \mathcal{W}_i$, and $\Delta_i$ is isometric to the preimage of a half-plane in an irreducible copy of $\mathbb{H}^2$ under the closest-point projection.
\end{itemize} 
\end{lemma}

For a quicker proof of this lemma, we remark that the natural symmetries of a triple in general position (given by Lemma \ref{lemma:duality}) preserve the walls.

\begin{lemma}
\label{lemma:dualwall}
Let $(F_i) \in \mathcal{F}^{(3)}$ be a positive triple of flags with walls $\mathcal{W}_i$ defined as before.
Then
\begin{enumerate}
    \item The order three symmetry $S$ permuting the flags verifies that $S \mathcal{W}_i = \mathcal{W}_{i+1}$ for $i=1,2,3$.
    \item The bilinear form $B_i$ inducing a duality fixing $F_i$ and permutting $F_{i-1}$ and $F_{i+1}$ verifies that $I_{B_i}$ fixes $\mathcal{W}_i$ and permutes $\mathcal{W}_{i-1}$ and $\mathcal{W}_{i+1}$.
\end{enumerate}
\end{lemma}

\begin{proof}
$(1)$ This follows from Lemma \ref{lemma:naturality}, as $S$ sends the data $(F_{i-1},F_{i+1},c_i)$ defining $\mathcal{W}_i$ to $(F_i, F_{i+2}, c_{i+1})$ defining $\mathcal{W}_{i+1}$.\\
$(2)$ By Lemma \ref{rmk:composition-dualities}, $D_{B_i}$ sends the spacelike geodesic corresponding to $\mathcal{W}_i$, to a geodesic asymptotic to $F_{i-1}$ and $F_{i+1}$.
Moreover, a computation using the expression of Lemma \ref{lemma:duality} shows that points in the component $c_i$ are sent again to points corresponding to $c_i$, thus, by the compatibility proven in Lemma \ref{rmk:composition-dualities} $(2)$, we get that $W_i$ is mapped to $W_i$.
Similarly, the spacelike geodesic corresponding to $W_{i+1}$ goes to a geodesic asymptotic to the same flags as $\mathcal{W}_{i-1}$.
There are two such walls, $\mathcal{W}_{i-1}$ and another one we will call $\mathcal{W}_{i-1}'$.
$\mathcal{W}_{i+1}$ cannot be mapped to $\mathcal{W}_{i-1}'$, as this wall contains parallel sets in timelike position with respect to parallel sets in $\mathcal{W}_i$.
\end{proof}

\begin{proof}[Proof of Lemma \ref{lemma:pingpongregions}]
By Lemma \ref{lemma:walls-separate}, the complement of $\mathcal{W}_i$ has two connected components, each isometric to the preimage of an irreducible half-plane via the closest-point projection.
Since the walls are disjoint, $\mathbb{X} \setminus \bigcup_{i=1}^3 \mathcal{W}_i$ has at most four components.
The order-three symmetry $S$ permuting the walls forces the existence of a component $\mathcal{D}$ bounded by the three walls.
\end{proof}

\begin{rmk}
\label{rmk:barycenter}
We can explicitly exhibit points in the domain $\mathcal{D}$.
For instance, the transformation $S$ permuting the three flags is a rotation, thus fixes a line $p_S \in \mathbb{RP}^2$, and leaves invariant a plane $\ell_S \in (\mathbb{RP}^2)^*$.
The parallel set $P_S$ defined by $(p_S,\ell_S) \in \mathbb{Y}$ is invariant under $S$, and thus belongs to the domain $\mathcal{D}$.
\end{rmk}

This leads us right into the definition of a parallel domain.

\begin{defi}
Given a positive triple of flags $(F_1, F_2, F_3)$, we denote the associated $\mathcal{D}$ of the lemma by $\mathcal{D}(F_1, F_2, F_3)$ and we call it \emph{parallel domain}. 
This parallel domain contains a parallel set $P_S$ invariant under the order-three symmetry $S$, which we call the \emph{barycenter} of the domain (as in Remark \ref{rmk:barycenter}).   
\end{defi}

\subsection{Distance estimates}
\label{subsec:estimates}

This section collects computations related to the distances between different parallel sets in our parallel domain, which will be crucial for what follows.
The core idea is that because the boundary of our domain is composed of parallel sets, and we can compute distances between parallel sets using cross ratios, we can obtain closed formulas.

To facilitate computations, we parametrize the geodesics defining the walls $\mathcal{W}_i$.
To do so, we make a (somewhat arbitrary) choice of basepoint.

\begin{figure}[ht]
    \centering
    \includegraphics[scale=.10]{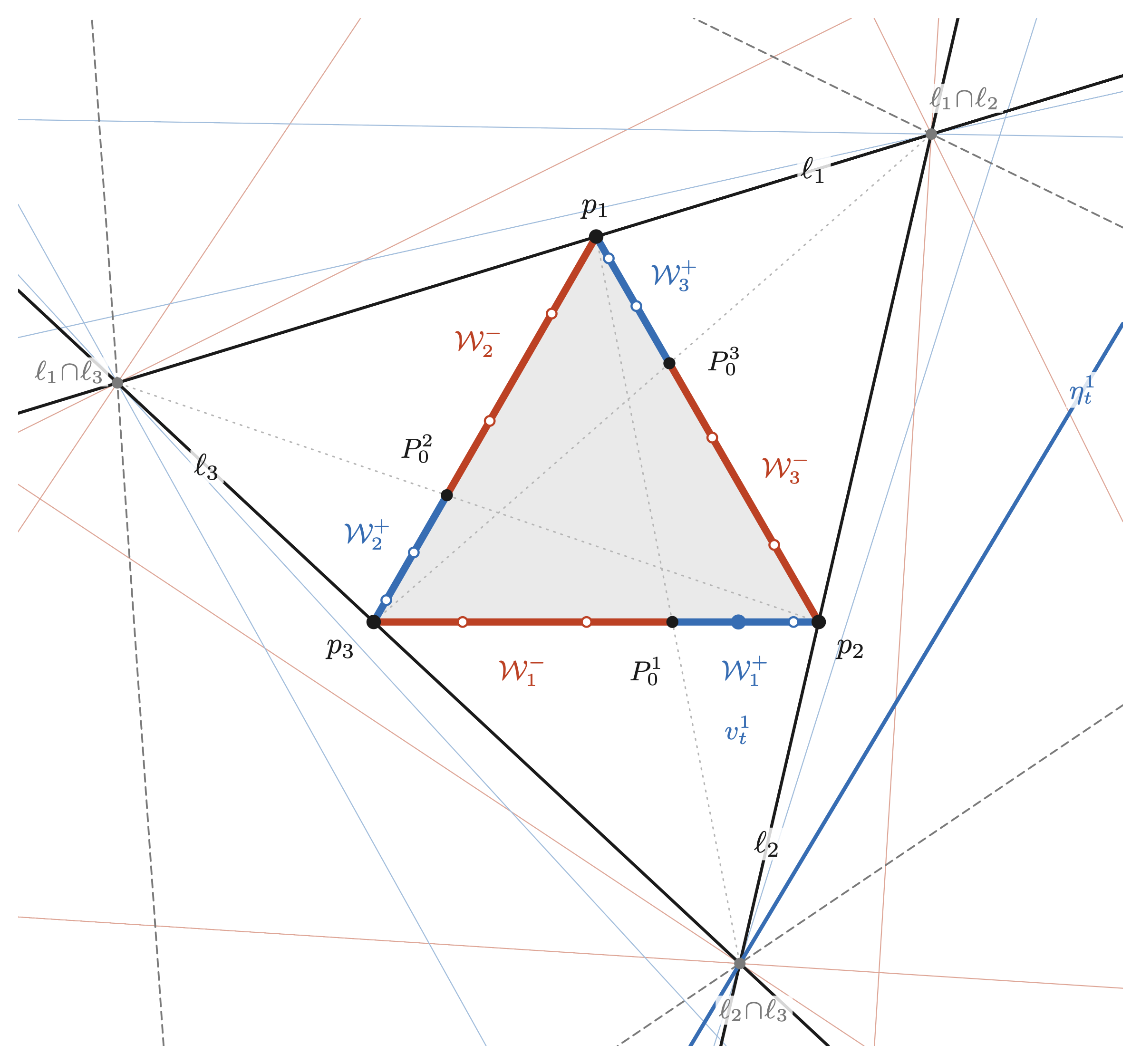}
    \caption{Parametrization of the walls $\mathcal{W}_i$.}
    \label{fig:separation}
\end{figure}

\begin{defi}
Let $(F_i) \in \mathcal{F}^{(3)}$ be a positive triple with parallel domain $\mathcal{D}$ bounded by walls $\mathcal{W}_i$, $i=1,2,3$.
\begin{itemize}
    \item For each $i \in \{1,2,3\}$, define the \emph{base parallel set} $P^i_0 \subset \mathcal{W}_i$ as the parallel set of the spacelike geodesic defining $\mathcal{W}_i$ whose point corresponds to $(p_{i-1} \oplus p_{i+1}) \cap (p_i \oplus (\ell_{i-1} \cap \ell_{i+1}))$ (see Figure \ref{fig:separation}).    
    \item We denote by $\mathcal{W}_i^+$ and $\mathcal{W}_i^-$ the connected components of $\mathcal{W}_i \setminus P_i$ corresponding to spacelike segments converging to $F_{i+1}$ and $F_{i-1}$ respectively.
\end{itemize}
\end{defi}

Using this decomposition, we will parametrize each parallel set of $\mathcal{W}_i$ by declaring that $P \in \mathcal{W}_i^{\pm}$ gets mapped to $\pm t$, for $ CR (P^i_0, P) = \cosh^2 (t)$.
By Lemma \ref{lemma:crossratiogeometry}, this equals the signed distance between $P$ and $P^i_0$ along the geodesic in the pseudo-Riemannian symmetric space.

\begin{prop}
\label{prop:estimates}
 Let $\mathcal{W}_i$, $i=1,2,3$ be the walls corresponding to the positive triple of flags with triple ratio $\tau \neq 0$.
 Then
    $$CR (P^i_t,P^{i+1}_{t'})
 =\frac{\left(1+e^{2t}+\tau e^{2t'+2t}\right)\left(1+\tau e^{2t}+\tau e^{2t'+2t}\right)}{4\tau e^{2t'+2t}}.$$
\end{prop}

\begin{proof}
Since $SL_3(\RR)$ acts by isometries on the symmetric space, it is enough to prove the proposition for a triple of flags of the form $F_i = ([e_i], [\alpha_i])$ for $\alpha_i =e^*_{i+1} + s e^*_{i-1}$ (a parametrization dual from the one in Lemma \ref{lemma:canonicalpos}).
In this case, we get $\tau^{-1/3} = s$.

In these coordinates, the point $\overline{p_{i-1} p_{i+1}} \cap (p_i \oplus (\ell_{i-1} \cap \ell_{i+1}))$ is defined by the vector $v^i_0 = e_{i+1} + s^2 e_{i-1}$. 
Proposition \ref{prop:asymptoticgeodesics} shows that:
$$v_t^i = e_{i+1} + e^{-2t} s^2 e_{i-1},\: \eta_t^i = \alpha_{i+1} + s e^{-2t} \alpha_{i-1}$$
define the point and lines corresponding to the parallel set $P_t^i$.

Direct computation shows that
$$\eta^i_t(v^{i+1}_{t'})=1+s^3 e^{-2t'}+s^3 e^{-2t'-2t},\: \eta^{i+1}_{t'} (v^i_t)=s \left( 1+e^{-2t'}+s^3 e^{-2t'-2t}\right),$$
while $\eta_t^i (v^i_t) = 2s^2 e^{-2t}$ and $\eta_{t'}^{i+1} (v^{i+1}_{t'}) = 2s^2 e^{-2t'}$.
Using that $s^3 = \tau^{-1}$ and the definition of the cross ratio, we get the desired expression.
\end{proof}

Particularly, since $$CR (P^i_t, P^{i+1}_{t'}) = \cosh^2 \left(d_\mathbb{X} (P^i_t, P^{i+1}_{t'})\right),$$ we get an expression for the distance between the parallel sets at the boundary of our domain.

We can turn this estimate into a qualitative statement useful for proving Theorem \ref{thmintro1}.

\begin{coro}
\label{coro:minimaldistance}
Let $(F_1,F_2,F_3)$ be a positive triple of flags inducing the walls $\mathcal{W}_i$ as before, we get that
\begin{equation}
\label{eq:minboundary}
\inf_{t' \in \mathbb{R}} \cosh \left(d_\mathbb{X} (P^i_t, P^{i+1}_{t'})\right) = \frac{\sqrt{1+e^{2t}} + \sqrt{1 + \tau e^{2t}}}{2},
\end{equation}
and the equality is attained at $t' = \frac{1}{4} \log \left( \frac{(1+e^{2t})(1+\tau e^{2t})}{\tau^2} \right) - t$.
Moreover, we get the following useful facts:
\begin{enumerate}
\item $d_\mathbb{X} (P_{t_n}^i,P_{t'_n}^{i+1}) \to 0$ if and only if $t_n \to -\infty$, $t'_n \to \infty$ and $t_n + t'_n \to -\frac{1}{2}\log \tau$.
\item There exists $\delta$ depending only on $\tau$ such that, for $i \neq j$,
$$d_\mathbb{X} (P^i_t,P^j_s) < \delta  \Rightarrow 
\begin{cases}
 j = i+1 \text{ and }t<0 ,s >0,\\
 j = i-1 \text{ and }t>0 ,s <0
\end{cases}.$$
\end{enumerate}
\end{coro}
\begin{proof}
Write $X = 1 + e^{2t}$, $Y = 1 + \tau e^{2t}$ and $Z = \tau e^{2t'+2t}$, then we can rewrite
$$CR(P_t^i, P_{t'}^{i+1}) = \frac{(X+Z)(Y+Z)}{4Z},$$
For fixed $t$, the function $f(Z)= (X/Z+1)(Y+Z) = X + Y + XY/Z + Z$ is strictly convex on $Z>0$, with minimum at $Z_m = \sqrt{XY}$. 
Plugging back this minimizer in terms of $t'$ and $t$, we obtain Equation \ref{eq:minboundary}.
Item $(1)$ follows immediately from this formula.

In what remains, we will prove item $(2)$.
It is enough to prove that there is a $\delta$ depending only on $\tau$ such that $d_\mathbb{X} (\mathcal{W}_i^+,\mathcal{W}_{i+1}),\: d_\mathbb{X} (\mathcal{W}_i^-,\mathcal{W}_{i+1}^-)> \delta$.
Note that the proof of Lemma \ref{lemma:dualwall}, shows that the duality $D_i$ exchanges $\mathcal{W}_i^+$ with $\mathcal{W}_i^-$. Therefore, these are the only cases to consider.

When $t \geq 0$, we get that
$$\inf_{t' \in \mathbb{R}} \cosh (d_\mathbb{X} (P_t^i,P_{t'}^{i+1})) \geq \frac{\sqrt{2}+ \sqrt{1 + \tau}}{2},$$
therefore, it suffices to work on the case $d_\mathbb{X} (\mathcal{W}_i^-,\mathcal{W}_{i+1}^-)$.

On the quadrant $\{ (t,t') : t \leq 0,\: t' \leq 0\}$, the function $CR (P^i_t,P^{i+1}_{t'})$ is proper, as $Z \to 0$ (as in the first paragraph) forces the cross ratio to go to infinity.
Thus, an infimum is attained.
This infimum cannot be attained in the interior of the quadrant, as if this were the case in $(t_0,t'_0)$ for $t_0,t'_0<0$, then this pair is realizing 
$$\inf_{t' \in \mathbb{R}} \cosh \left( d_\mathbb{X} (P^i_{t_0},P^{i+1}_{t'})\right)= \frac{\sqrt{1+e^{2t_0}} + \sqrt{1 + \tau e^{2t_0}}}{2}.$$
Since it is a global minimum, $\partial_t CR (P^i_t, P^{i+1}_{t'})$ vanishes at $(t_0,t'_0)$, but by the expression above, this is not the case.

Thus, the infimum is attained either at $t'=0$ or $t=0$.
Solving the minimization problem gives:
$$\cosh\left( d_\mathbb{X}(\mathcal{W}_i^-,\mathcal{W}_{i+1}^-) \right)
=\begin{cases}
\sqrt{\frac{(2+\tau)(1+2\tau)}{4\tau}}, & \tau\le\tau_0,
\:\text{when } (t,t')=(0,0),\\
\frac{\sqrt{2}+\sqrt{1+\tau^{-1}}}{2}, & \tau\ge\tau_0,
\:\text{when } \left(-\frac{1}{4}\log(2\tau(1+\tau)),\,0\right),
\end{cases}$$
where $\tau_0 = (\sqrt{3} - 1)/2$
In any case, this constant depends solely on $\tau$.
\end{proof}

The last estimate we will need concerns the distance from the barycenter $P_S$ of the domains to the walls.

\begin{lemma}
\label{lemma:distancebarycenter}
Let $P_t$ a parallel set in the wall $\mathcal{W}_i$, then
$$CR (P_S, P_t) = \frac{2}{3} \cosh (\log \tau^{1/6}) \left(\cosh (\log \tau^{1/6}) + \cosh (2t + \frac{1}{2} \log \tau)\right).$$
Particularly, $d_\mathbb{X} (P_S,P_t)$ is minimized when $t = -\frac{1}{4} \log \tau$.
\end{lemma}

\begin{proof}
This follows from a direct computation in the coordinates used to prove Proposition \ref{prop:estimates}, and from noticing that $P_S$ is associated with the vector $(1,1,1)$ and its orthogonal functional (with respect to the standard inner product).
\end{proof}

Notice that when $\tau =1$, this recovers that the inradius of an ideal triangle in $\mathbb{H}^2$ is $\cosh^{-1} (2/\sqrt{3}) = \frac{1}{2} \log 3$ by specializing to $\tau = 1$ and $t = 0$.
More generally, notice that the minimal distance from $P_S$ to the wall is
$$\inf_{t \in \mathbb{R}} \cosh (d_\mathbb{X} (P_S,P_t))= \frac{2}{\sqrt{3}} \cosh^{1/2} (\log \tau^{1/6}) \cosh (\log \tau^{1/12}),$$
which goes to infinity when $\tau \to 0,\infty$.
This gives intuition for the decay of the critical exponents we study later.

This lemma also suggests that our choice of center for the wall geodesic in the previous section was arbitrary.
It is better to shift the center of the geodesics by $\frac{1}{4} \log \tau$, since $P_{-\frac{1}{4}\log \tau}^i$ more clearly represents the foot of the altitude.
With this in mind, we can rewrite Proposition \ref{prop:estimates} as follows.

\begin{rmk}
\label{rmk:monotonicity}
Looking at the expression for $CR (P_t^i,P_{t'}^{i+1})$ given by Proposition \ref{prop:estimates} and defining $h = 2t + \frac{1}{2}\log \tau$ and $h' = 2t' + \frac{1}{2}\log \tau$, we can rewrite it as:
$$
\begin{aligned}
&\frac{\left(1+e^{2t}+\tau e^{2t'+2t}\right)\left(1+\tau e^{2t}+\tau e^{2t'+2t}\right)}{4\tau e^{2t'+2t}} =\\
 &\frac{\left(1 + e^{h+h'} + \tau^{-1/2} e^{h}\right)\left(1 + e^{h+h'} + \tau^{1/2} e^{h}\right)}{4e^{h+h'}},
\end{aligned}$$
One can verify that, for fixed $h$ and $h'$, the right-hand side is a monotonically increasing function of $\tau$ for $\tau \in [1,\infty)$.
\end{rmk}

%% file: sections/domains.tex
\section{Representations of $W_3$ and fundamental domains}
\label{domains:sec:representations}

This section has two goals.
First, we show that representations of the group $W_3$ into $SL_3(\RR)$ whose length-two elements are regular unipotent are in bijection with triples of flags which are pairwise in general position. 
We then show that the parallel domain $\mathcal{D}(F_1,F_2,F_3)$ defines a fundamental domain for the action of the representation $\rho : W_3 \rightarrow SL_3 (\mathbb{R})$ associated to the positive triplet $(F_1,F_2,F_3)$ in the symmetric space.

\subsection{The $SL_3 (\mathbb{R})$ character variety of $W_3$}
\label{subsec:representations}

Recall that $W_3 \leq \operatorname{Isom} (\mathbb{H}^2)$ is the group generated by reflections along an ideal triangle in the hyperbolic plane. 
A presentation for this group is
$$W_3 = \left\langle s_1, s_2, s_3 | s_i^2 \right\rangle$$
Notice that a representation of $W_3$ into $SL_3(\RR)$ is uniquely determined by a triple of involutions.

Starting from a triple of flags $(F_i) \in \mathcal{F}^{3}$ for which any pair is in general position, one can construct a representation of the ideal triangle reflection group $W_3$ into $SL_3(\RR)$. 
Such a triple give rise to three points $(\ell_{i+1} \cap \ell_{i-1}, p_{i+1} \oplus p_{i-1}) \in \mathbb{Y}$ (well defined by the pairwise transverse hypothesis). 
This corresponds to three involutions $\sigma_i$ via the correspondence described in Subsection \ref{prelim:involutions}. 

\begin{figure}[ht]
    \centering
    \includegraphics[width=0.4\textwidth]{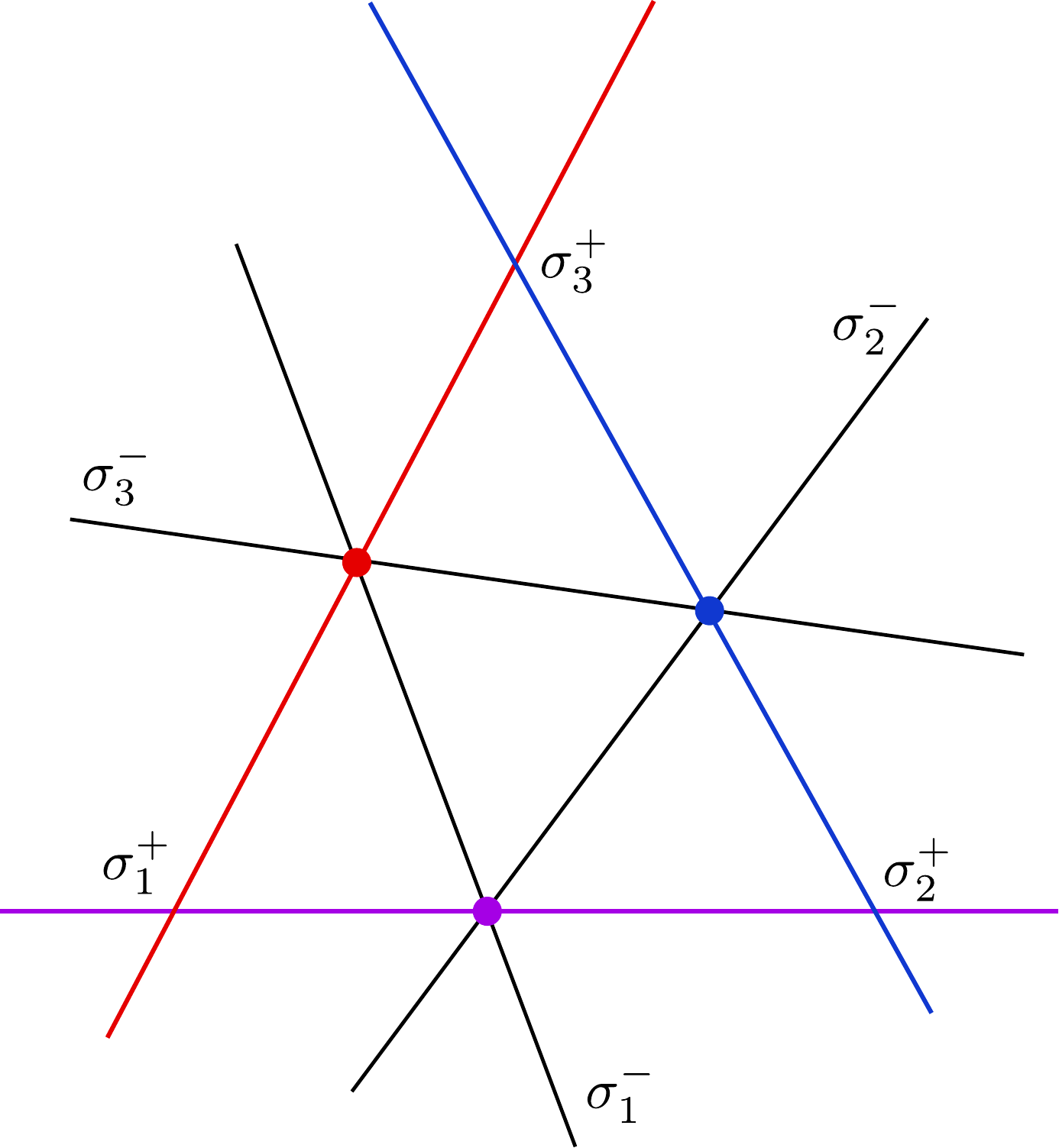}
    \caption{How to construct representations of $W_3$ from a flag. In purple we have $F_3$, in red $F_2$ and in blue $F_1$, one can see then how the three involutions are constructed.}
    \label{domains:fig:representation}
\end{figure}

These three involution define a representation $\rho_{(F_1, F_2, F_3)}: W_3 \rightarrow SL_3 (\RR)$ via $\rho_{(F_1, F_2, F_3)} (s_i) = \sigma_i$ (see Figure \ref{domains:fig:representation}).

By the way we constructed these representations, they verify that the composition $\sigma_{i-1} \sigma_{i+1}$ is a regular unipotent element fixing $F_i$.
In fact, one can observe the following:

\begin{lemma}
\label{lemma:unipotents}
Let $\sigma_1,\sigma_2 \in SL_3 (\mathbb{R})$ be a pair of involutions associated to $(p_i,\ell_i) \in \mathbb{Y}$, $i=1,2$.
Then, $\sigma_1 \sigma_2$ is conjugated to the regular unipotent element
$$u = \begin{pmatrix}
    1 & 1 & 0\\
    0 & 1 & 1\\
    0 & 0 & 1
\end{pmatrix},$$
if and only if $p_1 \neq p_2$, and $\ell_1 \cap \ell_2$ is contained in the line $p_1 \oplus p_2$.
If this happens, then the flag fixed by $\sigma_1 \sigma_2$ is $(\ell_1 \cap \ell_2, \overline{p_1p_2})$.
\end{lemma}

\begin{proof}
Suppose $\sigma_1 \sigma_2$ is regular unipotent, then it fixes a unique full flag.
Since the transformation has a unique invariant point, we have $\ell_1 \neq \ell_2$, and the point is $\ell_1 \cap \ell_2$.
Looking at the action in $(\mathbb{R}^3)^*$, the same argument shows that the unique invariant line is $\overline{p_1 p_2}$.
Since the invariant point lies on the invariant line, we obtain the converse direction.
The other direction is a direct computation. 
\end{proof}
With this lemma in mind we define the \emph{relative component of $W_3$} as
$$\operatorname{Hom}_{\text{rel}} (W_3, SL_3 (\RR)) = \{ \rho \in \operatorname{Hom}(W_3,SL_3 (\mathbb{R})) : \rho (s_{i-1}s_{i+1})\sim u, \: i = 1,2,3\},$$
where $g \sim h$ means that the matrices are conjugated. Notice that $SL_3 (\mathbb{R})$ acts on the relative character variety by conjugation.
Lemma \ref{lemma:unipotents} allows us to parametrize $\operatorname{Hom}_{\text{rel}} (W_3,SL_3 (\mathbb{R}))$ with triples of flags in general position.

\begin{prop}
\label{prop:parametrization}
The map $\operatorname{Hom}_{\text{rel}} (W_3,SL_3 (\mathbb{R})) \to \mathcal{F}^3$ sending a representation $\rho$, to the triple of flags $(F_i)$, for $F_i$ the unique flag fixed by $\rho (s_{i-1}s_{i+1})$, is a homeomorphism onto the triples of flags that are pairwise in general position.
Moreover, the subset of irreducible representations $\operatorname{Hom}^{\text{irr}}_{\text{rel}}(W_3,SL_3(\RR))\subset\operatorname{Hom}_{\text{rel}}(W_3,SL_3(\RR))$ is bijective to the space of triples in general position $\mathcal{F}^{(3)}$.
\end{prop}

\begin{proof}
Starting with $\rho \in \operatorname{Hom}_{\text{rel}}(W_3,SL_3 (\mathbb{R}))$, Lemma \ref{lemma:unipotents} claims that we can recover the points and planes defining the involutions $\rho (s_i)$ from the triple of flags.
This shows that the assignment $\operatorname{Hom}_{\text{rel}} (W_3,SL_3 (\mathbb{R})) \to \mathcal{F}^3$ is injective, and that the image verifies that any two flags are in general position.
The procedure mentioned above describes an inverse of this map, and thus proves surjectivity.
Continuity follows because the flags fixed by a unipotent element vary continuously with the element.

To prove the second item, it suffices to check that a representation is irreducible if and only if the three points of the flags are not aligned, and the three lines are not concurrent.
Notice that if either situation occurs, the representation is reducible (by taking either the $2$-plane corresponding to the line or the line corresponding to the concurrent point).
Thus, it suffices to see that the representation is irreducible whenever the triple is in general position.

Suppose $S \subset \mathbb{R}^3$ is an invariant subspace invariant under $\sigma_i = \rho (s_i)$.
It suffices to show that $S$ cannot be an invariant line, because if $S$ were two-dimensional, one could take the induced representation on the dual space, which would then have an invariant line (note that the condition imposed on the triple of flags is self-dual).
If $\dim S = 1$, then $S \in \mathbb{RP}^2$ is fixed by $\sigma_i$ for all i, which means $S$ is either the point invariant under $\sigma_i$, or lies in the invariant line.

We first show that $S$ cannot be a point invariant under some $\sigma_i$.
If so, $S$ cannot be a fixed point of $\sigma_j$, since the three lines of the flags are not concurrent.
Similarly, $S$ cannot be in the fixed line of $\sigma_j$, as this would imply that the line passing through $p_{j-1}$ and $p_{j+1}$ would be collinear with the line of the flag passing through one of the points, which is impossible as pairs of flags are in general position.

A similar argument shows that $S$ cannot be on the line fixed by $\sigma_i$.
If this happened, $S$ would lie in the intersection of the three invariant lines of the involution, which is impossible as the three points of the flags are not aligned.
\end{proof}

By Lemma \ref{lemma:canonicalpos}, $SL_3 (\mathbb{R})$-orbits of triples of flags in general position are in bijection with $\mathbb{R}\setminus \{0,-1\}$ via the triple ratio. 
The above proposition implies that the \textit{irreducible piece of the relative character variety} 
$$\mathcal{X}_{\text{rel}}^{\text{irr}} (W_3, SL_3 (\RR)) = \operatorname{Hom}^{\text{irr}}_{\text{rel}} (W_3, SL_3 (\RR))//SL_3 (\RR),$$ 
is isomorphic to $\mathbb{R} \setminus \{0,-1\}$.

We will refer to the connected component of $\mathcal{X}_{\text{rel}}^{\text{irr}} (W_3, SL_3 (\RR))$ corresponding to positive triples of flags as the \emph{relative Hitchin component} of $W_3$. 
This will be our object of study for the rest of the paper. 
The negative components also seem to have interesting properties and remain somewhat unexplored (see Appendix \ref{sec:negative}).

A quick calculation using the explicit parametrization of Lemma \ref{lemma:canonicalpos} shows that a triple of flags $(F_1,F_2,F_3)$ has triple ratio one if and only if the three flags are tangent to a common conic defined by a bilinear $(1,2)$ form $B$. 
By the construction of the involutions, the group $\rho_{(F_1,F_2,F_3)} (W_3)$ is contained in $SO(B)$.
We will call representations arising from such triples \emph{Fuchsian}. 

\subsection{Fundamental domains in the relative Hitchin component}

One of the main theorems in this paper is the following:
\begin{thm}
\label{thm:main}
Let $(F_1,F_2,F_3)$ be a positive triple of flags and $\rho : W_3 \rightarrow SL_3 (\mathbb{R})$ be the representation associated to this triple. 
Then, the parallel domain $\mathcal{D}(F_1,F_2,F_3)$ is a fundamental domain for the action of $\rho$ in the symmetric space $\mathbb{X}$.
\end{thm}

Before we explain the proof sketch, we need to verify that the involution acts as intended.

\begin{lemma}
\label{lemma:walls-invariantes}
Let $(F_i) \in \mathcal{F}^{(2)}$ be a pair of flags in general position, $\mathcal{W} = \mathcal{W}(F_1,F_2,c)$ a wall defined by this pair, and $\sigma$ the involution associated to the pair $(\ell_1 \cap \ell_2, p_1\oplus p_2) \in \mathbb{Y}$.
Then, the involution $\sigma$ preserves each parallel set in $\mathcal{W}$, and exchanges the connected components $\mathbb{X} \setminus \mathcal{W}$.
\end{lemma}

\begin{proof}
Pick an ellipsoid defined by a $(1,2)$ bilinear form $B$ tangent to the flags defining an irreducible slice $\mathcal{H}_B$ such that $\mathcal{W} = \pi^{-1} (\gamma)$, for $\gamma \subset \mathcal{H}_B$ a geodesic asymptotic to both flags, and $\pi$ the closest point projection.
The involution $\sigma$ preserves the ellipsoid by construction, fixes pointwise the projective segment corresponding to the geodesic $\gamma$, and exchanges the connected components of the complement.
Taking preimages of the closest-point projection gives the desired results.
\end{proof}

To simplify notation, write $\mathcal{D} = \mathcal{D}(F_1,F_2,F_3)$.
By fundamental domain, we mean that $\mathcal{D}$ has the following properties:
\begin{itemize}
    \item $\mathbb{X} = \bigcup_{g \in W_3} \rho (g) \overline{\mathcal{D}}$.
    \item $\overline{\mathcal{D}} \cap \rho (g) \overline{\mathcal{D}} \neq \emptyset$ if and only if the intersection is some $\mathcal{W}_i$, and $g = s_i$.
\end{itemize}

The proof of Theorem \ref{thm:main} follows the strategy used to prove the Poincar\'e polyhedron theorem.
We first build a tessellated space equipped with a Riemannian metric that is locally isometric to the symmetric space, and then prove that this local isometry is, in fact, an isometry.

\begin{defi}
In the setup of Theorem \ref{thm:main}, we define the \emph{tesselated space associated to $\rho$} as
$$\tess (\rho) := W_3 \times \overline{\mathcal{D}}/\sim,$$
where we define the equivalence relation $(g, x)\sim (gs_i, \rho (s_i)x)$ for $x \in \mathcal{W}_i$.
\end{defi}

The tessellated space comes equipped with a $W_3$-action by isometries by taking $h [(g,x)] = [(hg,x)]$.
Moreover, this action admits a fundamental domain of the form $\{ [(id,x)]: x \in \overline{\mathcal{D}}\}$ that we will denote as $\widetilde{\mathcal{D}}$.

There is a natural \emph{developing map} $\dev : \tess (\rho) \rightarrow \mathbb{X}$ defined via
$$\dev ([(g,x)]) = \rho (g) x.$$
This is well-defined thanks to the equivalence relations defining $\tess (\rho)$.
By construction, $\dev$ is a $\rho-$equivariant map, and clearly $\dev (\widetilde{\mathcal{D}}) = \overline{\mathcal{D}}$.
Using the developing map, we can pull back the structure of a Riemannian manifold from $\mathbb{X}$.

\begin{lemma}
The space $\tess (\rho)$ is a smooth manifold and admits a unique Riemannian metric that makes $\dev$ a local isometry.
\end{lemma}

\begin{proof}
If $x \in \mathcal{D}$ is an interior point, then any equivalence class $[(g,x)]$ has a unique representative, and therefore $\dev$ maps it bijectively to a neighborhood of the symmetric space.
If $x \in \mathcal{W}_i$, then choosing $U$ a $\rho (s_i)$-invariant neighborhood disjoint from $\mathcal{W}_{i-1} \cup \mathcal{W}_{i+1}$ (which exists because $\rho (s_i)$ is an involution), allows us to define $\mathcal{U} \subset \tess (\rho)$ as
$$\mathcal{U} = \{ [(g,x)] : x \in U \cap \overline{\mathcal{D}} \} \cup \{ [(gs_i, x)] : x \in U \cap \overline{\mathcal{D}}\}.$$
One can verify that $\dev|_{\mathcal{U}} : \mathcal{U} \rightarrow \mathbb{X}$ is bijective onto $U$.
Since $\rho(s_i)$ is an isometry for the Riemannian metric on $\mathbb{X}$, $\tess(\rho)$ inherits a Riemannian metric from the symmetric space.
\end{proof}

We will conclude the proof of Theorem \ref{thm:main} in two steps.
First, we prove that $\dev$ is globally injective (Proposition \ref{prop:pingpong}).
Then, we prove that $\tess (\rho)$ is geodesically complete (Proposition \ref{prop:geocomplete}).
Note that a local isometry between geodesically complete Riemannian manifolds is a covering.
In this case, the fact that $\mathbb{X}$ is non-positively curved manifold and simply connected finishes the proof, since this would imply the map is in fact a global isometry and thus a fundamental domain for $\tess (\rho)$ maps to a fundamental domain in $\mathbb{X}$.

Before proceeding to the proof, we note that global injectivity is needed to prove geodesic completeness.

\subsection{Proof of Theorem \ref{thm:main}}
We start by proving the global injectivity of the developing map.

\begin{prop}
\label{prop:pingpong}
 The map $\dev: \tess (\rho) \rightarrow \mathbb{X}$ is a globally injective local isometry.
\end{prop}

\begin{proof}
By equivariance of the developing map, it suffices to verify that for any non-trivial $g \in W_3$, we have $\rho(g) \mathcal{D} \cap \mathcal{D} = \emptyset$.
We will do so using a ping-pong argument adapted to the reflection group.
More precisely, let $\Delta_1,\Delta_2,\Delta_3$ be complementary regions of $\bigcup_{i=1}^3 \mathcal{W}_i$ as in Lemma \ref{lemma:pingpongregions}.
Note that the walls are disjoint by Proposition \ref{prop:criterion_disjointness}, and by Lemma \ref{lemma:walls-separate}, they separate the symmetric space and $\rho(s_i)$ exchanges the two connected components of $\mathbb{X} / \mathcal{W}_i$ (one of them being $\Delta_i$). All this implies that
\begin{equation}
\label{eq:pingpong}
\rho (s_i) \Delta_j \subset \Delta_i,\: \forall i \neq j.    
\end{equation}

Given $g \in \mathcal{W}_3$, the reduced word representing $g$ corresponds to a finite sequence of digits $1,2,3$ with no identical digits appearing consecutively.
Equation \ref{eq:pingpong} and an induction argument on the word length show that if the reduced form of $g \in W_3$ is of the form $g = s_i g'$, then $\rho (g) \mathcal{D} \subset \Delta_i$, implying $\rho (g) \mathcal{D} \cap \mathcal{D} = \emptyset$.
\end{proof}

As mentioned before, in order to prove the main theorem is enough to show the following proposition:

\begin{prop}
\label{prop:geocomplete}
The space $\tess (\rho)$ is geodesically complete.
\end{prop}

Here is where the estimates of Section \ref{subsec:estimates} come into play.
Before proving the proposition, notice the following:

\begin{lemma}
\label{lemma:infintersections}
Let $\widetilde{\gamma} : [0, T) \rightarrow \tess (\rho)$ be a unit-speed geodesic whose future cannot be extended. 
Then, it intersects infinitely many translates of the fundamental domain $\widetilde{\mathcal{D}}$ for the $W_3$-action.
\end{lemma}

\begin{proof}
Suppose by contradiction that $\widetilde{\gamma}$ intersects only finitely many translates of $\widetilde{\mathcal{D}}$. 
Then $\widetilde{\gamma}$ is necessarily contained in $\bigcup_{|g|\leq N} g \widetilde{\mathcal{D}}$ for $N$ sufficiently large so that all walls intersecting $\widetilde{\gamma}$ are interior to this set.
Up to taking a bigger $N$, we get that the closure of $\dev (\widetilde{\gamma})$ is contained in the interior of the open set $\dev (\bigcup_{|g|\leq N} g \widetilde{\mathcal{D}})$.
Since $\gamma = \dev \circ \widetilde{\gamma}$ is a geodesic of $\mathbb{X}$, and $\mathbb{X}$ is a complete metric space, $\gamma$ can be extended. 
A small extension of $\gamma$ beyond $T$ will still be included in $\dev (\bigcup_{|g|\leq N} g \widetilde{\mathcal{D}})$.

Since $\dev$ is an isometry onto its image, $\widetilde{\gamma}$ extends, which is a contradiction.
\end{proof}

This reduces the problem to understanding how geodesics can intersect infinitely many walls in finite time. 
We will use Corollary \ref{coro:minimaldistance} to understand this, and the following technical lemma which explains how parabolic elements act on spacelike geodesics.

\begin{lemma}
\label{lemma:paraboliclimit}
Let $\rho$ be a positive representation given by a triple of flags $(F_i)$, and let $u = \rho(s_is_{i+1})$, the unipotent that fixes the flag $F_{i - 1} = (p_{i-1}, \ell_{i-1})$. 
Then given $P$ a parallel set defining the wall $\mathcal{W}_i$ or $\mathcal{W}_{i+1}$, we get that $u^n(P) \to F_{i-1}$, for $|n| \to \infty$ uniformly over parallel sets of these walls. 
\end{lemma}

\begin{proof}
    In order to prove this lemma, we will use the normalization of the proof of Proposition \ref{prop:estimates}, that is, the triple of flags of the form $F_i = ([e_i], [\alpha_i])$ with $\alpha_i = e_{i+1}^* + se_{i-1}^*$.
Without loss of generality, we will prove it for $i = 1$ (using the order three symmetry gives the proof for the other cases). 
From Proposition \ref{prop:asymptoticgeodesics}, we have that the elements of $\mathbb{Y}$ given by 
$$v_t^1 = e_{2} + e^{-2t} s^2 e_{3},\qquad \eta_t^1 = \alpha_{2} + s e^{-2t} \alpha_{3} = (s + se^{-2t})e_1^* + s^2e^{-2t}e_2^* + e_3^*$$
parametrize elements of $\mathcal{W}_1$.
Under the chosen normalization, a quick calculation shows that 
$$u^n = 
\begin{pmatrix}
-2n + 1 & -2ns & 0\\
2n/s & 2n + 1 & 0\\
\frac{2n^2 - 2n + 2n^2s^3}{s^2} & \frac{2n^2 + 2n^2s^3 + 2ns^3}{s} & 1
\end{pmatrix}.$$

Putting both things together, we get 
$$u^nv^1_t = -2nse_1 + (2n + 1)e_2 + (s^2e^{-2t} + \frac{2n^2 + 2n^2s^3 + 2ns^3}{s})e_3,$$
and, because $u^n$ stabilizes $\alpha_3$, it follows that
$$u^n \eta^1_t = se^{-2t}\alpha_3 + \frac{2n^2 + 2n + 2n^2s^3 + 2ns^3 + s^3}{s^2}e_1^* + \frac{2n^2 + 2n^2s^3}{s}e_2^* + e_3^*. $$

Moreover, because $s$ is positive, 
$$\left|s^2e^{-2t} + \frac{2n^2 + 2n^2s^3 + 2ns^3}{s}\right| \ge \frac{2n^2 + 2n^2s^3 + 2ns^3}{s}.$$
Note that the term containing $e_3$ grows at least quadratically with $n$ for every $t$, while the other terms grow linearly. This implies that $[u^nv^1_t]$ converges to $[e_3]$ uniformly on $t$. 
As for the $\eta_t^1$, note that
$$u^n\eta^1_t = \left(se^{-2t} + \frac{2n^2 + 2n^2s^3}{s^2} \right)\alpha_3 + \frac{2n + 2ns^3 + s^3}{s^2}e_1^*  + e_3^*.$$
Again, the fact that $s > 0$ implies that 
$$\left|se^{-2t} + \frac{2n^2 + 2n^2s^3}{s^2} \right| \ge \frac{2n^2 + 2n^2s^3}{s^2},$$
which in turn implies that $[u^n\eta^1_t]$ converges uniformly to $[\alpha_3]$, by an argument analogous to the one above.
\end{proof}

Now we are ready to prove the proposition, which in turn prove Theorem \ref{thm:main}.

\begin{proof}[Proof of Proposition \ref{prop:geocomplete}]
Suppose, by contradiction, that there exists a geodesic that cannot be extended $\widetilde{\gamma}: [0, T) \rightarrow \tess (\rho)$.
Observe that $\gamma = \dev \circ \widetilde{\gamma}$ can be extended in the symmetric space, and $\gamma (T) \in \partial \dev (\tess (\rho))$.
Note that $\gamma$ cannot be contained in a wall because of Lemma \ref{lemma:infintersections}.

Let $(t_k')$ be the increasing sequence of times for which $\gamma (t_k')$ intersects some translate of a wall $\mathcal{W}_{i_k'}$.
We will call the label $i_k'$ the \emph{index} of the wall.
By Remark \ref{rmk:finitelymanyint}, the sequence $t_k'$ can only accumulate at $T$.

The arc $\gamma$ may, a priori, intersect the same translate of a wall twice in a row. We define an alternative sequence $(t_k)$ as follows: set $t_1 = t_1'$ and, for $k>1$, we set $t_{k+1}$ as the smallest $t'_j$ bigger than $t_k$ such that $\gamma(t_k)$ and $\gamma(t'_j)$ belong to different walls.
Notice that Lemma \ref{lemma:infintersections} ensures that this sequence is infinite.
Moreover, Remark \ref{rmk:finitelymanyint} implies that $ t_k$ converges to $ T$. 
We will now introduce some notation: 
\begin{itemize}
\item $i_k$ is the index of the wall containing $\gamma (t_k)$.
\item $g_k \in W_3$ is the element for which $\gamma((t_{k+1} - \varepsilon, t_{k+1})) \subseteq \rho(g_k)\mathcal{D}$, for some $\varepsilon > 0$.
\item $P_{r_k}^{i_k}$ for $r_k \in \mathbb{R}$ is the parametrized parallel set of $\mathcal{W}_{i_k}$ whose translate contains $\gamma (t_k)$ (as in Section \ref{subsec:estimates}).
\end{itemize}
By definition of $t_{k+1}$, on the interval $(t_k,t_{k+1})$ the geodesic $\gamma$ meets no wall other than the one containing $\gamma(t_k)$, so it stays in the two tiles adjacent to that wall.
In particular, $\gamma(t_k)$ and $\gamma(t_{k+1})$ lie on the walls $\rho(g_k)\mathcal{W}_{i_k}$ and $\rho(g_k)\mathcal{W}_{i_{k+1}}$ of the same tile $\rho(g_k)\overline{\mathcal{D}}$.
By the definition of $\tess(\rho)$ and Proposition \ref{prop:pingpong}, the only tiles adjacent to a wall $\rho(g)\mathcal{W}_j$ are $\rho(g)\overline{\mathcal{D}}$ and $\rho(gs_j)\overline{\mathcal{D}}$.
Since both $\rho(g_k)\overline{\mathcal{D}}$ and $\rho(g_{k+1})\overline{\mathcal{D}}$ are adjacent to the wall containing $\gamma(t_{k+1})$, we get
\begin{equation}
\label{eq:consecutive-tiles}
g_{k+1} \in \{ g_k,\ g_k s_{i_{k+1}} \}.
\end{equation}
Moreover, every tile met by $\gamma$ after $t_1$ is $\rho(g_k)\overline{\mathcal{D}}$ or $\rho(g_k s_{i_k})\overline{\mathcal{D}}$ for some $k$, so Lemma \ref{lemma:infintersections} implies that $\{g_k\}$ is unbounded.

The second item of Corollary \ref{coro:minimaldistance} implies the existence of $\delta>0$ depending only on $\tau$ such that, if $i \neq j$,
$$d(P^i_t,P^j_{t'}) < \delta  \Rightarrow 
\begin{cases}
    j = i+1 \text{ and }t<0 ,t' >0,\\
    j = i-1 \text{ and }t>0 ,t' <0
\end{cases}.$$
Let $k_0$ be such that for any $k \geq k_0$, we get the bound $d_\mathbb{X} (\gamma (t_k),\gamma (t_{k+1}))  = |t_k - t_{k+1}| < \delta$.
Note that
$$d_\mathbb{X} (P_{r_k}^{i_k}, P_{r_{k+1}}^{i_{k+1}}) \leq d_\mathbb{X} (\gamma (t_k),\gamma (t_{k+1})) = |t_k - t_{k+1}|<\delta.$$
Therefore, our choice of $k$ implies that,
\begin{itemize}
    \item if $r_k < 0$, then $i_{k+1} = i_k + 1$ and $r_{k+1} > 0$.
    \item if $r_k > 0$, then $i_{k+1} = i_k -1$ and $r_{k+1} < 0$.
\end{itemize}
Using both conditions simultaneously, we get that 
$$i_{k_0 + j} = 
\begin{cases}
i_{k_0},\: j\text{ even}\\
i_{k_0} \pm 1,\: j\text{ odd}.
\end{cases}$$

The sign in the previous formula is determined by the sign of $r_{k_0}$.
We will assume that $r_{k_0} < 0$, the other case is analogous.
By \eqref{eq:consecutive-tiles} and the previous formula, for $k \geq k_0$ every $g_k$ is either of the form
$$g_k = g'(s_{i_{k_0}} s_{i_{k_0}+1})^{n_k} \quad \text{or} \quad g_k = g'(s_{i_{k_0}} s_{i_{k_0}+1})^{n_k} s_{i_{k_0}},$$
for some $n_k \in \mathbb{Z}$, where $g' = g_{k_0}$.
For $k = k_0 + 2j$ we have $i_k = i_{k_0}$, and since $\rho(s_{i_{k_0}})\mathcal{W}_{i_{k_0}} = \mathcal{W}_{i_{k_0}}$, in both cases
$$\gamma(t_k) \in \rho(g')\rho(s_{i_{k_0}} s_{i_{k_0}+1})^{n_k}\mathcal{W}_{i_{k_0}}.$$
Up to taking a subsequence, we may assume that $|n_k| \to \infty$.

We started the proof by assuming that $\gamma(T) \in \partial \dev(\tess(\rho))$. 
Up to composing with $\rho (g')^{-1}$, the previous discussion shows that there is a sequence of parallel sets $P_{r_k}^{i}$, for $i = i_{k_0}$ contained in $\mathcal{W}_i$, points $x_k \in P_{r_k}^i$ and a sequence $(n_k)$ going to infinity such that $\rho (s_i s_{i+ 1})^{n_k} (x_k) \to x \in \mathbb{X}$.
We will see that this produces a contradiction.

Let $(p_k, \ell_k) \in \mathbb{Y}$ be such that $P(p_k, \ell_k) = \rho (s_i s_{i+ 1})^{n_k} P_{r_k}^i$. 
Because of Lemma \ref{lemma:convergencia-paralelos}, $(p_k, \ell_k)$ has a convergent subsequence in $\mathbb{Y}$. 
On the other hand, Lemma \ref{lemma:paraboliclimit} implies that $(p_k, \ell_k)$ converges to a flag $F_{i-1} \in \mathcal{F}$, which is a contradiction.
\end{proof}

%% file: sections/criticalexponent.tex
\section{Application to critical exponents}
\label{sec:critical-exponent}

In this section, $\rho_\tau$ denotes the representation determined by a positive triple of flags $(F_1, F_2, F_3)$ with triple ratio $\tau$. 
We will see how the distance estimates obtained for the fundamental domains allow us to compute asymptotics for the critical exponent $\rho_\tau$ in the symmetric space.
Recall that this number is defined as
$$h_\mathbb{X} (\rho_\tau) = \limsup_{T\to \infty} \frac{1}{T} \log \# \{g \in W_3 : d_\mathbb{X} (o,\rho (g)o) \leq T\},$$
for $o \in \mathbb{X}$ some choice of basepoint (this number is independent of this choice).
We will show:

\begin{thm}
\label{thm:limite-exp}
With the setup from before, we get that
$$\lim_{\tau \to +\infty} h_\mathbb{X} (\rho_\tau) = \frac{1}{2}.$$
\end{thm}

We remark that $\frac{1}{2}$ is precisely the critical exponent of the group generated by a regular unipotent, thus $h_\mathbb{X} (\rho_\tau) \geq \frac{1}{2}$ for any $\tau > 0$.
We also prove something stronger: the critical exponent associated with the Hilbert length of $\rho_\tau$ tends to $\frac{1}{2}$ as $\tau \to \infty$ (see Theorem \ref{thm:criticalhilbert}).

\subsection{Chains of parallel sets}

Since the critical exponents of $\rho_\tau$ are always bounded below by $1/2$, our result follows from producing upper bounds for the critical exponents in terms of $\tau$.
Upper bounds can be obtained with lower bounds for $d(o,\rho_\tau (g) o)$ in terms of $\tau$.
To do this, we have the freedom of picking $o$ in the parallel set $P_S \subset \mathcal{D}$, and notice that $d_\mathbb{X} (P_S,\rho_\tau (g) P_S) \leq d (o,\rho_\tau (g) o)$.
In this section, we will focus on producing lower bounds for the left-hand side of this inequality.

Using Theorem \ref{thm:main}, we can control $d_\mathbb{X} (P_S,\rho_\tau (g)P_S)$ by looking at which parallel sets of walls a minimizing geodesic from $P_S$ to $\rho_\tau (g)P_S$ must cut.

\begin{defi}
Let $\rho_\tau : W_3 \rightarrow SL_3 (\mathbb{R})$ a representation as before, and $g \in W_3$ with reduced form $g = s_{j_1}\ldots s_{j_n}$.
\begin{itemize}
\item A \emph{chain} from $P_S$ to $\rho_\tau (g)P_S$ is a sequence $P_1,\ldots, P_{n}$, such that
$$P_k = \rho_\tau (g_{k-1}) P_{t_k}^{j_k},$$
where $g_k = s_{j_1}\ldots s_{j_k}$, and $P_{t_k}^{j_k}$ is a parallel set of the wall $\mathcal{W}_{j_k} \subset \partial \mathcal{D}$ parametrized as in Section \ref{subsec:estimates}.
Denote the set of such chains as $\mathcal{C}_\tau (g)$.
\item The \emph{chain translation length} $D_\tau (g)$ from $P_S$ to $\rho_\tau (g) P_S$ is defined as
\begin{equation}
\label{eq:chaintranslength}
D_\tau (g) = \inf_{(P_i) \in \mathcal{C}_\tau (g)}\sum_{i=0}^{n} d_\mathbb{X} (P_i,P_{i+1}), 
\end{equation}
where we declare that $P_0 = P_S$ and $P_{n+1} = \rho_\tau (g) P_S$.
A chain is \emph{minimizing} if it realizes the infimum. 
\end{itemize}
\end{defi}

Our motivation for defining $D_\tau (g)$ comes from the following lemma.

\begin{lemma}
\label{lem:comparacion-chain}
Given $g \in W_3$ and $\tau \in \mathbb{R}^+$, then 
$$D_\tau (g)\leq d_\mathbb{X} (P_S,\rho_\tau (g) P_S).$$
When $\tau = 1$, equality is achieved, and $D_1 (g) = d_{\mathbb{H}^2} (z_B, g z_B)$, where we think of $W_3 \leq \operatorname{Isom}(\mathbb{H}^2)$, and $z_B$ is the barycenter of the ideal triangle defining the reflection group.
\end{lemma}

\begin{proof}
Take a geodesic $\gamma$ realizing $d_\mathbb{X} (P_S, \rho_\tau (g)P_S)$.
By Remark \ref{rmk:finitelymanyint}, this geodesic intersects translates of the walls in $\mathcal{D}$ at times $0<d_1'< \ldots < d_N' < \ell (\gamma)$, and two consecutive times lie in a translate $\rho_\tau (g_k') \mathcal{D}$ of the parallel domain.
Notice that $(g_k')_{k=1,\ldots,N}$ defines a path in the Cayley graph of $W_3$ from the identity to $g$.
This Cayley graph is a tree of cycles, i.e., a graph obtained from an infinite ternary tree by doubling each edge.
In particular, geodesics from the origin are non-backtracking paths.

The geometry of this tree implies that if we denote by $g_k = s_{j_1}\ldots s_{j_k}$ obtained from truncating $g = s_{j_1}\ldots s_{j_n}$, then $(g'_j)$ necessarily passes through all $g_k$.
Moreover, when the path reaches $g_k$ for the first time, then it must necessarily come from $g_{k-1}$.

This allows us to inductively define times $0 < d_1 <\ldots < d_{n} < \ell (\gamma)$ such that $d_k = d'_{l_k}$, for $\gamma (d'_{l_k})$ the first time the segment between $\gamma (d'_{l_k})$ and $\gamma (d'_{l_k+1})$ is adjacent to the tile $\rho_\tau (g_k) \mathcal{D}$.
Since the tile transversed by $\gamma$ before $\gamma (d_k)$ corresponds to $g_{k-1}$, we get that $\gamma (d_k) \in P_k =\rho_\tau (g_{k-1}) P_{t_k}^{j_k}$ for any $k$.
In particular, $\{P_k\}$ forms a chain.
Thus
$$\sum_{i=0}^{n} d_\mathbb{X} (P_i,P_{i+1}) \leq \sum_{i=0}^{n} d(\gamma (d_i), \gamma (d_{i+1})) = d_\mathbb{X} (P_S,\rho_\tau (g) P_S).$$
Taking the infimum gives the desired inequality.

For the second part, note that when the triple ratio is one, all the parallel sets are orthogonal to the same irreducible slice (as in Section \ref{subsec:irredparallel}).
In particular, the distances between the parallel sets coincide with the distances of their respective projections to the slice.
\end{proof}

Note that any term in the sum in Equation \eqref{eq:chaintranslength}, except the first and last, can be computed using Proposition \ref{prop:estimates}, and the other two can be computed using Lemma \ref{lemma:distancebarycenter}.
In particular, the left-hand side depends explicitly on the triple ratio and the times $(t_j)$ defining a chain.

Notice that a sequence of times $(t_i)_{i=1,\ldots,n}$ determines a unique chain from $P_S$ to $\rho_\tau (g) P_S$.
The estimates from Section \ref{subsec:estimates} allow us to construct a bijection between chains $\mathcal{C}_\tau (g)$ and $\mathcal{C}_{\tau'} (g)$ in a way that varies nicely with $\tau$.
More precisely, define:

\begin{equation}
\label{eq:bijection-pattern}
c_{\tau,\tau'} : \mathcal{C}_\tau (g) \rightarrow \mathcal{C}_{\tau'} (g),
\end{equation}
mapping a chain between $P_S$ and $\rho_\tau (g) P_S$ corresponding to times $(t_i)$, to a chain from $P_S$ to $\rho_{\tau'} (g) P_S$ associated times $(t_i')$, where we declare that
$$t_i' = t_i + \frac{1}{4} \log \frac{\tau}{\tau'}, \forall i \in \{1,\ldots,n\}.$$
This choice of times is motivated by Remark \ref{rmk:monotonicity}, which implies the following.

\begin{rmk}
\label{rmk:ctautau}
Given $(P_i) \in \mathcal{C}_\tau (g)$, we get that $(P_i') = c_{\tau,\tau'} (P_i)$ is a sequence verifying that $d(P_i', P_{i+1}') < d(P_i,P_{i+1})$ whenever $\tau' < \tau$ (and the reverse inequality if $\tau' > \tau$), for $\tau, \tau' \ge 1$.
\end{rmk}

This leads us to the following monotonicity result:

\begin{prop}
\label{prop:monotone}
Given $g \in W_3$, then the function $\tau \to D_\tau (g)$ is monotinically increasing in $[1,\infty)$.
More precisely, we get the inequality
$$D_\tau (g) - D_{\tau'} (g) \geq \log \frac{\cosh (\log \tau^{1/6})}{\cosh (\log \tau'^{1/6})},$$
for any $1 \leq \tau' < \tau$.
\end{prop}

\begin{proof}
Let $\{P_k\}$ be a minimizing chain for $P_S$ and $\rho_{\tau} (g) P_S$, and let $c_{\tau,\tau'} (P_i) = (P_k')$ the corresponding chain between $P_S$ and $\rho_{\tau'} (g) P_S$, as constructed in the previous section. By setting $P_0 = P_S$ and $P_{n + 1} = \rho_\tau(g)(P_S)$, we get
$$D_\tau (g) - D_{\tau'} (g) \geq \sum_{k=0}^n d(P_k,P_{k+1}) - d(P'_k,P'_{k+1}).$$
Remark \ref{rmk:monotonicity} claims that every term, except the first and last, is strictly positive.
Thus, it suffices to bound the first and last terms.

Notice that for our choice of $t_1$ and $t_1'$, Lemma \ref{lemma:distancebarycenter} implies
$$\frac{CR(P_S,P_1)}{CR(P_S,P_1')} \geq \frac{\cosh (\log \tau^{1/6})}{\cosh (\log \tau'^{1/6})}.$$
Since 
$$\frac{CR(P_S,P_1)}{CR(P_S,P_1')} \leq e^{2 (d_\mathbb{X} (P_S,P_1) - d_\mathbb{X} (P_S,P_1'))},$$
taking logs and repeating the argument for the other end finishes the proof.
\end{proof}

This suggests that critical exponents in the symmetric space, or, as we will later see, the Hilbert metric, should decrease monotonically in $\tau$.
However, the additive constant between $D_\tau (g)$ and $D_{\tau'} (g)$ is not enough to obtain this result.

\subsection{The Farey lamination and bounds in terms of the Coxeter length}
\label{subsec:farey}

We can think of $W_3$ as the subgroup of $\text{Isom} (\mathbb{H}^2)$ generated by reflections along edges of the ideal triangle $\Delta_0$ with vertices at $0$,$1$,$\infty$. Label the edges of $\Delta_0$ by $e_1, e_2, e_3$ according to which generator of $W_3$ fixes them.
The orbit of these edges defines the \emph{Farey lamination}.

Observe that the equivariance of the walls and the existence of the fundamental domain (Theorem \ref{thm:main}) give a $W_3$-equivariant bijection between the orbit of the walls by the action of the group and the Farey lamination in $\HH^2$, by sending $ge_i \to \rho_\tau(g)\mathcal{W}_i$.

We will use the notion of \emph{cutting sequences}, introduced by Series in \cite{series1985modular}, to study the combinatorics of elements in the group through fundamental domains. We adapt her definition to match our setting.

\begin{figure}[ht!]
\includegraphics[width=0.8\linewidth]{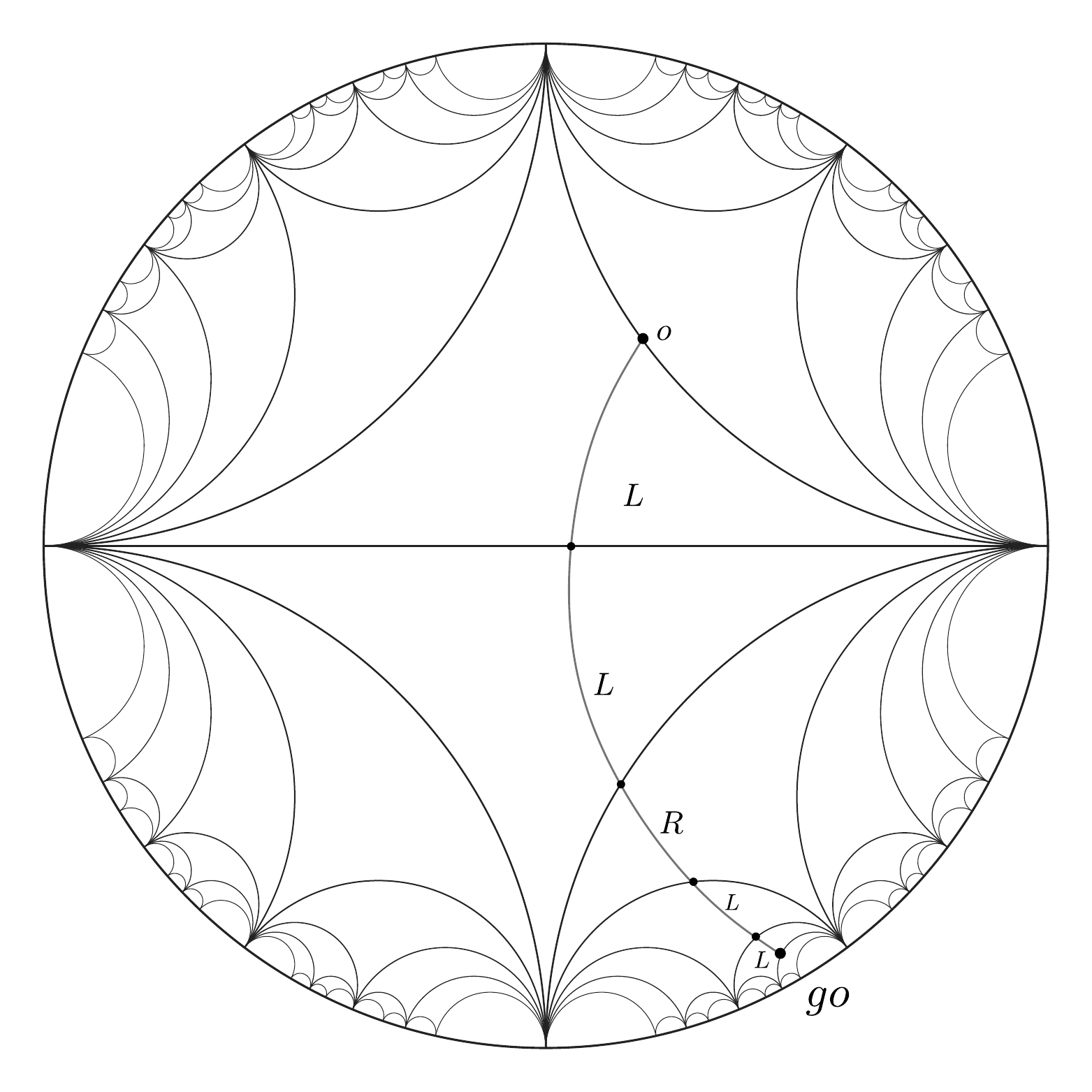}
\caption{An explanation on how a cutting sequence starting from a point $o$ is constructed.}
\label{fig:farey} 
\end{figure}

Suppose $\Delta_0$ is the ideal triangle with vertices $0$, $1$ and $\infty$ in the upper half-plane model. 
Let $i \in \mathbb{H}^2$, $g \in W_3$ and take $\gamma$ be the oriented geodesic in $\mathbb{H}^2$ from $i$ and $g i$.
Observe that this geodesic intersects different leaves of the lamination at points $\gamma (0), \gamma(t_1), \ldots, \gamma(t_n)$, where $\gamma(t_n) = g i$.
For each arc from $\gamma (t_i)$ to $\gamma (t_{i+1})$, the two leaves of the lamination containing $\gamma(t_i)$ and $\gamma(t_{i+1})$ intersect at an ideal point in $\partial \HH^2$.
We associate the symbol $L$ or $R$ with the point according to whether it lies to the left or right of $\gamma$. 
This allows us to associate to $g$, a \emph{cutting sequence} $L^{k_1} R^{k_2}\ldots L^{k_j}$, where each $k_i \in \mathbb{N}^+$, except possibly $k_1$ and $k_j$ (that may be zero) (see Figure \ref{fig:farey}).

\begin{defi}
Given $g \in W_3$ with cutting sequence $L^{k_1} R^{k_2} \ldots L^{k_j}$, we will call each $k_i$ the \emph{depth of a cusp excursion}.
We will call the number $C(g)$ of alternations between powers of $L$ and powers of $R$ the \emph{Coxeter length} of $g$. A visualization of a cutting sequence can be seen in Figure \ref{fig:farey} 
\end{defi}

\begin{rmk}
\label{rmk:turning}
Consider the Farey pattern generated by the ideal triangle joining the vertex $0$, $1$ and $\infty$ of $\partial \HH^2$. The stabilizer of this pattern in $\mathbb{H}^2$ is, up to index two, $SL_2 (\mathbb{Z})$, and up to index 6, $W_3$.
Moreover, the only elements of this stabilizer leaving $\{\Re z = 0\}$ invariant are generated by the reflection along the geodesic and the rotation of angle $\pi$. This group has order 4.
\end{rmk}

Setting 
$$L = 
\begin{pmatrix}
1 & 1 \\
0 & 1    
\end{pmatrix},
\:
R = 
\begin{pmatrix}
1 & 0 \\
1 & 1
\end{pmatrix}$$
Associates a cutting sequence with an element of $PSL_2 (\mathbb{Z})$.
These elements allow us to estimate translation lengths. The following lemma is classical, but we include its proof for completeness:

\begin{lemma}
\label{lem:cutting-sequence}
Let $g \in W_3$ with reduced form $s_{j_1}\ldots s_{j_n}$ and cutting sequence $L^{k_1}R^{k_2}\ldots L^{k_j}$, corresponding to $M = L^{k_1}R^{k_2}\ldots L^{k_j} \in PSL_2(\mathbb{Z})$.
If $\Re(gi) > 0$, then $gi = Mi$, and more generally, $d(i, gi) = d(i, Mi)$.
\end{lemma}

\begin{proof}
The rotation $r(z) = -1/z$ fixes $i$, preserves the Farey pattern, and preserves orientation, so the segment from $i$ to $r g(i)$ is the $r$-image of the one from $i$ to $gi$ and has the same cutting sequence. 
Hence we may assume $\Re(gi) > 0$.

Let $\gamma$ be the geodesic joining $i$ and $g i$, and let $0= t_0 < t_1 <\ldots < t_n = d(i, gi)$ be the sequence of times at which the geodesic intersects the Farey pattern in a leaf $\gamma (t_k) \in e_k$, and denote by $M_k$ the product of the first $k$ left and right turns $L, R$ determined by the segment from $\gamma (0)$ to $\gamma (t_k)$. 
Notice that $e_0$ is the geodesic $\{\Re z = 0\}$.
We will prove by induction that $M_k e_0 = e_k$ for all $k$.

The case $k=0$ is trivial.
Assume $e_j = M_j e_0$ for any $j \leq k$ and denote by $\Delta_j$ the triangle bounded by $e_j$ and $e_{j+1}$.
$\Delta_j = M_j \Delta_0$ for any $0 \leq j \leq k-1$ by inductive hypothesis. 
Since $M_{k-1}$ preserves the tessellation and the orientation, and maps $e_0$ to $e_{k-1}$ and $\Delta_0$ to $\Delta_{k-1}$, it maps the left (resp. right) edge of $\Delta_0$ seen from $e_0$ to the left (resp. right) edge of $\Delta_{k-1}$ seen from $e_k$, and likewise for the triangles beyond them. 
This shows that $M_{k-1}^{-1} \Delta_k$ agrees with $L$ or $R$ applied to $\Delta_0$, depending on what letter corresponds to the transition from $e_k$ to $e_{k+1}$, which shows the desired result.

In particular $gi \in e_n = M e_0$, so $h = M^{-1} g$ maps $i$ into $e_0$. 
Since $h$ preserves the tessellation and distinct edges have disjoint interiors, $h e_0 = e_0$. 
An element of the stabilizer of the Farey tessellation that preserves $e_0$ must fix $i$, so in particular $gi = Mi$.
\end{proof}

Now we will estimate $D_\tau (g)$ in terms of the Coxeter length $C(g)$.
Corollary \ref{coro:minimaldistance} bounds the minimal distance between $P^i_t$ and an adjacent wall in terms of the function 

$$F_\tau (x) = \cosh^{-1} \left(\frac{\sqrt{1+e^{2x}}+\sqrt{1+\tau e^{2x}}}{2}\right).$$
We get the following elementary property from this function:
\begin{lemma}
Given $\tau >0$, the function $F_\tau$ is a strictly convex function.    
\end{lemma}

And thus, we get the following corollary:

\begin{coro}
Let $P_\cdot^i$ denote a parallel set in $\mathcal{W}_i$ at the boundary of our fundamental domain.
Then
$$d(P_s^i, P_t^{i+1}) + d(P_t^{i+1},\rho_\tau (s_{i+1}) P_r^{i-1}) \geq 2 \cosh^{-1} \left(\frac{\sqrt{1+\tau^{-1/2}} + \sqrt{1 + \tau^{1/2}}}{2}\right),$$
for any $s,t,r \in \mathbb{R}$. 
The same inequality holds for $P_s^i$, $P_t^{i-1}$ and $P_t^{i+1}$.
\end{coro}

\begin{proof}
We know that $d(P_t^{i+1},\rho_\tau (s_{i+1}) P_r^{i-1}) \geq F_\tau (t)$ by Corollary \ref{coro:minimaldistance}.
Aplying a duality, we get that $d(P_s^i, P_t^{i+1}) \geq F_\tau (-t - \frac{\log \tau}{2})$.
Now, by the convexity of the function $F_\tau$, we get
$$d(P_s^i, P_t^{i+1}) + d(P_t^{i+1},\rho_\tau (s_{i+1}) P_r^{i-1}) \geq F_\tau (t) + F_\tau (-t-\log \tau) \geq 2 F_\tau \left(-\frac{\log \tau}{4}\right),$$
which is the result we wanted.
By symmetry, we obtain a similar result for $i$, $i-1$, and $i+1$.
\end{proof}

This corollary controls the length of geodesics with one endpoint in $\mathcal{W}_{i + 1}$ (resp. $\mathcal{W}_{i - 1}$) and another in $\sigma_i\mathcal{W}_{i-1}$ (resp. $\sigma_i\mathcal{W}_{i+1}$). 
Whenever we get a transition from $R$ to $L$ (or vice versa) in the cutting sequence of a word $g \in W_3$, the situation from the previous sentence translates.
This implies the following.

\begin{thm}
\label{thm:coxeteriality}
Let $\tau >0$, and $g \in W_3$, then
$$D_\tau (g) \geq (C(g) - 4)\cosh^{-1} \left(\frac{\sqrt{1+\tau^{-1/2}} + \sqrt{1 + \tau^{1/2}}}{2}\right).$$
\end{thm}

\begin{proof}
Let $g_k = s_{i_1}\ldots s_{i_k}$ be the truncation of $g = s_{i_1}\ldots s_{i_n}$ of length $k$.
The claim is that $C(g)$ is equal to the number of $k$ such that $g_{k+3} = s_{i \pm 1} s_i s_{i \mp 1} g_k$ up to an additive error of four.
This follows automatically by composing with $g_k^{-1}$ and looking at the cutting sequence for the geodesic joining the barycenter $z_B$ and $s_{i \pm 1} s_i s_{i\mp 1} z_B$ (notice that we are maybe forgetting about two alternations at the beginning or at the end).
This shows that the parallel sets in the chain $P_{k}$, $P_{k+1}$, and $P_{k+2}$ are translates of parallel sets to which the previous corollary applies.
Since each term $d(P_k,P_{k+1})$ belongs to at most two sequences from $k$ to $k+3$, we get the desired inequality.
\end{proof}

This estimate is good for words whose length is proportional to its Coxeter length.
However, it is a poor estimate for words that make long excursions along cusps.
We balance this by using Theorem \ref{thm:coxeteriality} and Proposition \ref{prop:monotone} simultaneously.

\subsection{Bound on the excursions}

This section gives a lower bound for the translation length in the Fuchsian representation in terms of the cusp excursions, which Proposition \ref{prop:monotone} then uses to bound the chain translation length.

\begin{lemma}
\label{lem:distance-estimate}
Let $g \in W_3$ with cutting sequence $L^{n_1} R^{n_2} \ldots L^{n_k}$, then
$$d (i , gi) \geq 2 \sum_{j=1}^{k} \log \max\{n_j, 1\}.$$
\end{lemma}

\begin{proof}
Because of Lemma \ref{lem:cutting-sequence} we have that
$d(i, gi) = d(i,L^{n_1}R^{n_2}\ldots L^{n_k} i)$.
    
A simple argument by induction shows that at least one entry of the matrix $L^{n_1}R^{n_2}\ldots L^{n_k}$ is greater than or equal to $\max \{n_1,1 \} \ldots \max\{n_k,1\}$ (notice that $n_i \geq 1$ for all $i$ except maybe the first and last ones).
This shows that
$$\|g\|_2 \geq  \prod_{j=1}^k \max\{n_j,1\},$$
where the norm is taken with respect to the standard inner product, as $SO(2)$ stabilizes $i$. 
Notice that $\|M\|_2 \leq e^{d(i, Mi)/2} = e^{d(i, gi)/2}$. 
Taking logs gives the desired result. 
\end{proof}

We also remind the reader of the following classical computation.

\begin{lemma}
\label{lemma:exponente-parabolico}
 Let $A =
    \begin{pmatrix}
        1 & 1\\
        0 &1
    \end{pmatrix}$, then the critical exponent $h_{\langle A\rangle} = \frac{1}{2}$.
\end{lemma}

\begin{proof}
Notice that $\cosh d(i, A^n i) = 1 + n^2/2$. 
As $n \to \infty$, we get that $d (i, A^n i) \sim 2\log n$. 
Thus, for every $\delta > 1$ we have that for $n$ big enough, 
$$2\log n \le d(i, A^n i) \le 2\log (\delta n).$$
This implies that the critical exponent is $1/2$.
\end{proof}

The following proposition will be crucial.
Informally, it quantifies how many elements with high Coxeter index are in a hyperbolic ball of radius $R$.

\begin{prop}
\label{prop:manhatan}
 Let $\zeta (x) = \sum_{n \geq 1} \frac{1}{n^x}$ denote the Riemann zeta function, then for any $y > \frac{1}{2}$ and $x > \log \zeta (2 y)$ we get that $\sum_{g \in W_3} e^{-x C(g) - y d(i,gi)}$ converges.
\end{prop}

\begin{proof}
Recall from the proof of Lemma \ref{lem:cutting-sequence} that an element in $W_3$ coincides with the element related to its cutting sequence, up to an element of the stabilizer of the Farey pattern that fixes $(0, \infty)$. 
The stabilizer of $(0, \infty)$ has order 4. This, along with Lemma \ref{lem:distance-estimate}, implies that we have the following chain of inequalities:
    \begin{align*}
    \sum_{g \in W_3} e^{-x C(g) - y d(i,gi)} &= \sum_{r \geq 0} e^{-x r} \sum_{C(g) = r} e^{-y d(i,gi)}\\
    &\leq 4 \sum_{r \geq 0} e^{-x r} \left( \sum_{n_1,\ldots,n_{r+1}} (\max\{n_1, 1\}n_2\ldots n_r\max\{n_{r+1}, 1\})^{-2y} \right).
    \end{align*}
 Notice that $$\sum_{n_1,\ldots,n_{r+1}} (\max\{n_1, 1\}n_2\ldots n_r\max\{n_{r+1}, 1\})^{-2y} \leq C\zeta(2y)^{r+1}$$
 for some $C > 0$.
 Putting everything together, we obtain a geometric series that converges whenever the proposition's conditions are met.
\end{proof}

\subsection{Proof of Theorem \ref{thm:limite-exp}}

We now compute the asymptotics of the critical exponent in the symmetric space.
The key result is the following. 

\begin{prop}
\label{prop:laposta}
When $\tau$ goes to $\infty$, the critical exponent $H(\tau)$ of the Dirichlet series $\sum_{g \in W_3} e^{-s D_\tau (g)}$ defined as
$$H(\tau) = \inf \{s >0 :  \sum_{g \in W_3} e^{-s D_\tau (g)}< \infty \},$$
converges to $\frac{1}{2}$.
\end{prop}

\begin{proof}
A calculation analogous to Lemma \ref{lemma:exponente-parabolico} shows that $D_\tau((s_1s_2)^n) \sim 2\log n$, which implies that $H(\tau) \ge \frac{1}{2}$.

By the last paragraph, to show convergence it suffices to check that for any $s > 1/2$, there is a $\tau(s)$ such that the Dirichlet series converges for any $\tau > \tau(s)$.

We have two bounds for $D_\tau (g)$:
\begin{itemize}
\item Theorem \ref{thm:coxeteriality} implies that $D_\tau (g) \geq C(g) G(\tau)$, for $G(\tau) = F_\tau(-\frac{1}{4}\log \tau)$, which goes to infinity when $\tau \to \infty$.
\item Proposition \ref{prop:monotone} implies that $D_\tau (g) > D_1 (g)$, and $D_1 (g) = d_{\mathbb{H}^2} (i,gi)$.
\end{itemize}
Particularly, we get that $D_\tau (g) \geq t d_{\mathbb{H}^2} (i,gi) + (1-t) C(g) G(\tau)$ for any convex combination $t \in (0,1)$.
    
Take $s > \frac{1}{2}$, and consider $t$ such that $st > \frac{1}{2}$ and $t < 1$.
We get that:
$$\sum_{g \in W_3} e^{-s D_\tau (g)} \leq \sum_{g \in W_3} e^{-s (1-t) G(\tau) C(g) - s t d(i,gi)}.$$
Since $G(\tau) \to \infty$ when $\tau \to \infty$, we obtain that $s (1-t) G(\tau) > \log \zeta (2 st)$ for sufficiently large $\tau$.
Proposition \ref{prop:manhatan} tells us that the right-hand side converges. 
Thus, the series on the left converges, and $s$ is greater than the critical exponent for sufficiently large $\tau$.
Since this holds for any $s > \frac{1}{2}$, we get the desired claim.
\end{proof}

Theorem \ref{thm:limite-exp} is a direct consequence.

\begin{proof}[Proof of Theorem \ref{thm:limite-exp}]
Notice that $h_\mathbb{X} (\rho_\tau)$ is the critical exponent of the Dirichlet series:
$$\sum_{g \in W_3} e^{-s d_\mathbb{X} (o,\rho_\tau (g)o)}.$$
The presence of unipotents implies that $h_\mathbb{X} (\rho_\tau) \geq \frac{1}{2}$.
On the other hand, $d_\mathbb{X} (o,\rho_\tau (g)o) \geq D_\tau (g)$, and therefore
$$1/2 \leq h_\mathbb{X} (\rho_\tau) \leq H(\tau).$$ 
The previous proposition implies $\lim_{\tau \to \infty} H(\tau) = \frac{1}{2}$, which finishes the proof.
\end{proof}

\subsection{Asymptotics for the Hilbert critical exponent}

In this section,  we prove that the critical exponents of the representations $\rho_\tau$ with respect to the Hilbert functional converge to $1/2$ as $\tau \to \infty$. This is stronger than Theorem \ref{thm:limite-exp}.

We denote by $\mathfrak{a}$ the Cartan subspace defined in Section \ref{subsec:regular}, $o \in \mathbb{X}$ a preferred choice of basepoint in the symmetric space, and $\kappa : SL_3 (\mathbb{R}) \rightarrow \mathfrak{a}^+$ the Cartan projection induced by $o$ as in Section \ref{subsec:finslercartan}.

We will refer to $\alpha^H \in \mathfrak{a}^*$ defined via:
$$\alpha^H (x,y,z) = \frac{1}{2} (x-z),$$
as the \emph{Hilbert functional}.
This functional allows us to define the \emph{Hilbert critical exponent} of $\rho_\tau$ as
$$h_H (\rho_\tau) = \limsup_{T\to \infty} \frac{1}{T} \log \# \{g \in W_3: \alpha^H (\kappa (\rho_\tau (g))) \leq T\}.$$

This number has two interpretations.
Geometrically, this quantity is the growth rate of orbits with respect to a Finsler norm on the symmetric space defined by the functional $\alpha^H$, as introduced in \cite{kapovichleeb2018}.
Dynamically, this quantity coincides with the topological entropy of a reparametrization of the geodesic flow on the orbifold associated to $W_3$ (see \cite{kimoh25,blayac24}).
In this situation, this reparametrization has a nice geometric realization in terms of the Hilbert geometry of a convex domain associated to $\rho_\tau$, thanks to Vinberg's theorem \cite{vinberg1971discrete} and Benoist's work \cite{benoist01}.

We now show that Proposition \ref{prop:laposta} implies the following.

\begin{thm}
\label{thm:criticalhilbert}
Let $h_H(\rho_\tau)$ be the Hilbert critical exponent. 
Then 
$$\lim_{\tau \to \infty} h_H(\rho_\tau) = 1/2.$$
\end{thm}

This is stronger than Theorem \ref{thm:limite-exp} because $\alpha (\kappa (\rho_\tau (g))) \leq d_\mathbb{X} (o,\rho_\tau (g))$.
The inequality follows from Cauchy-Schwarz and the fact that $\| \kappa (\rho_\tau (g))\| = d_\mathbb{X} (o,\rho_\tau (g)o)$.

The only ingredient we need is the following lemma.

\begin{lemma}
\label{lemma:P_shilbertlength}
Suppose the basepoint $o$ belong to the barycenter $P_S \subset \mathcal{D}$, then 
$$d_\mathbb{X} (P_S, \rho (g)P_S) \leq \alpha^H (\kappa (g)) + \log 4.$$
\end{lemma}

Theorem \ref{thm:criticalhilbert} comes as a direct consequence of this Lemma.

\begin{proof}[Proof of Theorem \ref{thm:criticalhilbert}]
By Lemma \ref{lem:comparacion-chain}, we get that $D_\tau (g) \leq \alpha^H (\kappa) + \log 4$.
This bound implies that $h_H (\rho_\tau) \leq H(\tau)$.
Since $h_H (\rho_\tau) \geq 1/2$ for any $\tau$ due to the presence of a unipotent element, we get the desired asymptotic.
\end{proof}

In order to prove the lemma, we start by noticing the following intermediary steps that will help us compute the distance from $P_S$ to $\rho_\tau (g) P_S$:

\begin{lemma}
Let $P_S$ denote the barycenter of the parallel domain associated to $\rho_\tau$. 
Then $P_S$ and $\rho_\tau (P_S)$ are in spacelike position for any $\tau > 0$.
\end{lemma}

\begin{proof}
Consider the set
$$I = \{\tau \in \mathbb{R}^+: P_S \text{ and }\rho_\tau (g) P_S\text{ are in spacelike position}\}.$$
This set is non-empty as $P_S$ and $\rho_1 (P_S)$ are in spacelike position. 
Since spacelike configurations are open in $\mathbb{Y}^2$, the set $I$ is open.
It suffices to see closedness.

Suppose $(\tau_n) \in I^\mathbb{N}$ is a sequence converging to $\tau \in \mathbb{R}^+$.
If $\tau \not\in I$, then the parallel sets $P_S$ and $\rho_\tau P_S$ cannot be in a non-degenerate configuration of $\mathbb{Y}^2$, as these are open.
Therefore, $d_\mathbb{X} (P_S,\rho_\tau P_S) = 0$, i.e., the parallel sets are either asymptotic or they cross.
This is impossible by Lemma \ref{lemma:distancebarycenter}.
\end{proof}

Thus, we can use the formula $CR (P_S,\rho_\tau (P_S)) = \cosh^2 (d_\mathbb{X} (P_S,\rho_\tau (P_S)))$ to estimate the distance.

\begin{proof}[Proof of Lemma \ref{lemma:P_shilbertlength}]
Let $([v_S],[\alpha_S]) \in \mathbb{Y}$ associated to $P_S$.
Suppose $v_S$ is unitary for $q$. Because, by definition, $v_S \perp^q \ker \alpha_S$, up to rescaling $\alpha_S$ we may assume that $\alpha_S(v) = q(v, v_S)$. Thus, we have that
$$CR(P_S,\rho (g)P_S) = q(g^{-1} v_S ,v_S) q( g v_S, v_S ) \leq \|g^{-1}\|_q\|g\|_q,$$
where we used the Cauchy-Schwarz inequality.
Note that $\|g\|$, and $\|g^{-1}\|^{-1}$ are the first and last entries of $\kappa (g)$. 
Taking logarithms and using the previous remark gives the desired inequality.
\end{proof}

%% file: sections/negative.tex
\section{The negative component}
\label{sec:negative}

This paper concerns the study of one component of the relative character variety $\mathcal{X}_{\text{rel}}^{\text{irr}} (W_3,SL_3 (\mathbb{R}))$. 
Namely, the one corresponding to positive triples of flags.
In this section, we argue that some representations in the other components exhibit an interesting behaviour that differs from the studied component.
Note that the two negative components are identified through any duality, as a duality sends triple ratio $\tau$ to $\tau^{-1}$.

Consider the normalization of the $SL_3(\RR)$-orbit of a triple of flags in general position given by
$$F_i (s) = ([e_{i+1} + s e_{i-1}], [e_i^*]), \: i =1,2,3.$$
A computation shows that $TR(F_1, F_2, F_3) = s^3$. 
Notice that the transformation $S \in SL_3(\RR)$ defined by $S e_i = e_{i+1}$ verifies that $S F_i (s) = F_{i+1} (s)$.

Following the procedure of Section \ref{subsec:representations}, each triple $(F_1 (s),F_2 (s),F_3 (s))$ defines a representation $\rho_s : W_3 \rightarrow SL_3 (\mathbb{R})$ defined via:

$$\rho_s (s_1) = 
\begin{pmatrix}
1 & -2/s & -2s\\
0 & -1 & 0\\
0 & 0 & -1
\end{pmatrix},
\rho_s (s_2) = S \rho_s (s_1) S^{-1},\: \rho_s (s_3) = S^{-1} \rho_s (s_1) S.$$
In this appendix, we focus on the case $s < 0$.
We start by observing that this component cannot preserve a convex domain, as happens for the positive component. 
More generally, the representation cannot be relatively Anosov, a generalization to higher rank of geometrically finite representations (see \cite{zhu2025relatively}). 
This is shown in the following Lemma, as every non-parabolic, infinite order element in a relative Anosov representation needs to be loxodromic.
\begin{lemma}
\label{lemma:coxeter}
When $s < 0$ and $ s \neq-1$, the Coxeter element $s_3 s_2 s_1$ has a complex eigenvalue.
Particularly, $\rho_s$ cannot be relatively Anosov with respect to the peripheral system of the ideal reflection orbifold.
\end{lemma}

\begin{proof}
Recall that a cubic polynomial of the form
$$x^3 + ax^2 + bx + c,$$
has a complex root if and only if its discriminant $\Delta = 18 abc -4a^3c + a^2 b^2 -4b^3 -27c^2$ is less than zero.
Computing the characteristic polynomial using the explicit matrices from above, we get that
$$\chi (\rho_s (s_1s_2s_3)) = X^3 + (9+ \frac{8}{s^3}) X^2 - (8s^3 + 9)X - 1.$$
Therefore:
$$\Delta (s) = 18 (9+\frac{8}{s^3})(8s^3+9) + 4 (9+\frac{8}{s^3})^3 + (9 + \frac{8}{s^3})^2(9+8s^3)^2+4(8s^3+9)^3-27.$$
A simple calculation shows that $\Delta(s) \leq 0$ for all $s < 0$ and equals $0$ only at $s = -1$.
\end{proof} 

There is one distinguished boundary case of a negative triple of flags, namely, when 
$$TR(F_1,F_2,F_3) = -1.$$
In this situation, Menelaus theorem claims that the three points of the flags are aligned, or the three lines are concurrent at a point.
The associated representation is therefore reducible, and after some computations one notices that it is not faithful.
This is not very auspicious for the negative component.

However, notice that when $s = -2$ (or $s=-1/2$), the representation $\rho_{-2}$ takes values on $SL_3 (\mathbb{Z})$.
Particularly, $\rho_{-2} (W_3) = \Gamma$ is automatically discrete.
Moreover, one can check that the group $\Gamma$ projects surjectively to $SL_3 (\mathbb{Z}/5 \mathbb{Z})$, and therefore, $\Gamma$ is Zariski dense in $SL_3 (\mathbb{R})$ (see for instance \cite{Lubotzky1999generation}).
Thus, we ask the following question:

\smallskip
\textbf{Question:} Is $\rho_{-2} : W_3 \rightarrow SL_3 (\mathbb{Z})$ faithful?
\smallskip

Since free groups do not have property (T), a positive answer to the previous question would show that $\Gamma$ is thin (i.e., Zariski dense and of infinite index).

If $\Gamma$ were thin, it would satisfy two interesting properties:
first, since the Coxeter element maps to a complex-loxodromic element, \cite[Theorem 1.4]{dey2025remarks} implies that $\Gamma$ acts minimally on the full flag manifold $\mathcal{F}$.
To our knowledge, no finitely generated, infinite-volume discrete subgroups of $SL_3 (\mathbb{R})$ are known to act minimally on the full flag manifold. 
Second, this would answer a question by Kapovich \cite[Problem 3.16]{kapovich2023list}, who asked whether every thin subgroup of $SL_3 (\mathbb{Z})$ is (relatively) Anosov.

We have extensive computational evidence pointing towards an affirmative answer to the previous question.
First, no element of length $\leq 45$ with respect to the standard generators of $W_3$ is mapped to the identity.
Second, we observed that the only elements mapped to unipotents are those in the peripheral system, i.e., those conjugated to $s_i s_j$.
The remaining elements map to real or complex loxodromic elements.

\begin{figure}[ht]
    \centering
    \includegraphics[scale=.8]{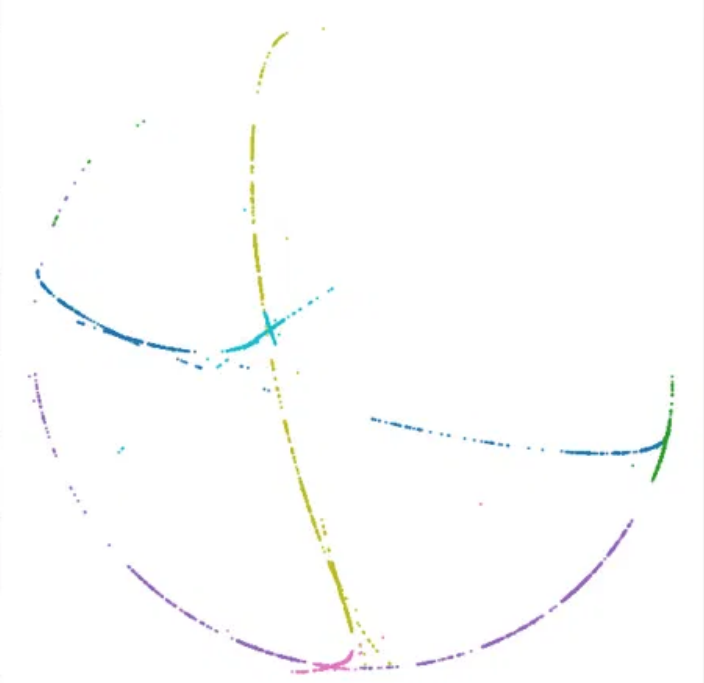}
    \caption{Each color corresponds to the orbit of a vertex of a fundamental domain of $W_3$, with mixed colors corresponding to mixed accumulation. We use a stereographic model of $\mathbb{RP}^2$ for the picture.}
\label{fig:negative}
\end{figure}

We can use the structure of the Farey tessellation to construct a pattern of flags associated to $\Gamma$.
Namely, the endpoints of geodesics of the Farey lamination are in bijection with $\mathbb{P}^1 (\mathbb{Q})$.
Any point in $\mathbb{P}^1 (\mathbb{Q})$ is the fixed point of a parabolic element in $W_3 \leq \text{Isom} (\mathbb{H}^2)$. 
We can construct an assignment $\mathbb{P}^1 (\mathbb{Q}) \rightarrow \mathcal{F}$ by choosing a unipotent $u$ fixing a given point in $\mathbb{P}^1 (\mathbb{Q})$ and sending it to the flag fixed by $\rho_{-2} (u)$.
Notice that the images are in fact composed of rational flags.
Figure \ref{fig:negative} depicts this assignment.

We do not expect to extend this map continuously to $\mathbb{P}^1$ to obtain a continuous Cannon-Thurston map $\xi: \mathbb{P}^1 \rightarrow \mathcal{F}$ because of the existence of a complex loxodromic element.

